\documentclass[11pt,reqno,a4paper]{amsart}
\usepackage[margin=2.0cm,top=2cm,bottom=2cm]{geometry}
\usepackage{setspace}
\usepackage{amsmath,amssymb}
\usepackage{xcolor} 
\usepackage{amsthm}
\usepackage{amsfonts}
\usepackage{subcaption}
\usepackage{graphicx}
\usepackage{csquotes}
\usepackage{enumerate}
\usepackage{soul}

\usepackage[colorlinks=true]{hyperref}
\hypersetup{urlcolor=blue, citecolor=red}
\makeatletter
\newcommand*{\rom}[1]{\expandafter\@slowromancap\romannumeral #1@}
\makeatother
\newcommand{\N}{\mathbb{N}}
\newcommand{\R}{\mathbb{R}}

\newcommand{\ve}{\varepsilon}
\newcommand{\bp}{\begin{pmatrix}}
\newcommand{\ep}{\end{pmatrix}}
\newcommand{\intt}{\int_{\T^d}}
\newcommand{\p}{\partial}
\newcommand{\iinttr}{\iint_{\T^d\times\R^d}}
\newcommand{\ov}{\overline}

\newtheorem{theorem}{Theorem}[section]
\newtheorem{lemma}{Lemma}[section]
\newtheorem{corollary}[lemma]{Corollary}

\newtheorem{proposition}[lemma]{Proposition}

\theoremstyle{definition}
\newtheorem{definition}[lemma]{Definition}

\theoremstyle{remark}

\newtheorem{remark}[lemma]{Remark}

\numberwithin{equation}{section}

\DeclareMathOperator*{\argmin}{arg\,min}

\newcommand{\T} {\mathbb T}
\newcommand{\pa}{\partial}

\newcommand{\lt}{\left}
\newcommand{\rt}{\right}
\newcommand{\bq}{\begin{equation}}
\newcommand{\eq}{\end{equation}}

\newcommand{\bfb}{{\bf b}}
\newcommand{\bfc}{{\bf c}}

\newcommand{\bfX}{{\bf X}}

\newcommand{\bbN}{\mathbb N}

\newcommand{\bbS}{\mathbb S}

\newcommand{\calD}{\mathcal D}
\newcommand{\calE}{\mathcal E}

\newcommand{\calJ}{\mathcal J}

\newcommand{\calL}{\mathcal L}
\newcommand{\calM}{\mathcal M}

\newcommand{\calP}{\mathcal P}

\newcommand{\calR}{\mathcal R}
\newcommand{\calS}{\mathcal S}

\newcommand{\weakto}{\rightharpoonup}

\newcommand{\intr}{\int_{\R^d}}

\newcommand{\dx}{\textnormal{d}x}
\newcommand{\dv}{\textnormal{d}v}

\newcommand{\dt}{\textnormal{d}t}

\usepackage{soul}
\usepackage[normalem]{ulem}
\usepackage{cancel}

\newcommand{\dd}{\mathrm{d}}

\newcommand{\AC}{\operatorname{AC}}
\newcommand{\lefthalfcup}{\mathbin{\vrule height 1.6ex depth 0pt width
0.13ex\vrule height 0.13ex depth 0pt width 1.3ex}}

\usepackage{accents}
\newcommand{\ddt}{\frac{\textnormal{d}}{\dt}}

\renewcommand{\sout}[1]{}
\renewcommand{\cancel}[1]{}

\begin{document}

\title[Lagrangian solutions to the ionic Vlasov--Poisson system]{Global existence of Lagrangian solutions to the ionic Vlasov--Poisson system}

\author{Young-Pil Choi}
\address{Department of Mathematics, Yonsei University, Seoul 03722, Republic of Korea}
\email{ypchoi@yonsei.ac.kr}

\author{Dowan Koo}
\address{Mathematical Institute, University of Oxford, Oxford OX2 6GG, United Kingdom}
\email{dowan.koo@maths.ox.ac.uk}

\author{Sihyun Song}
\address{Department of Mathematics, Yonsei University, Seoul 03722, Republic of Korea}
\email{ssong@yonsei.ac.kr}

\date{\today}



\begin{abstract} 
In this paper, we study the Cauchy problem of the ionic Vlasov--Poisson (IVP) system for a broad class of initial data. To this end, we establish well-posedness and stability estimates for the Poisson--Boltzmann equation $-\Delta \Phi = \rho -e^\Phi$ provided $\rho \in L_+^p$ with $p>1$, whereas previous approaches required $p>\frac{d}{2}$. Building upon these results, we construct global-in-time renormalized/Lagrangian solutions to the IVP system.
\end{abstract}


\maketitle

\setcounter{tocdepth}{1}
%
%
%
%
%

\section{Introduction}
We study the ionic Vlasov--Poisson (IVP) system, which is also known as the Vlasov--Poisson system with massless electrons:
\begin{align}
\label{iVP}
\begin{cases}
    \p_t f_t + v \cdot \nabla_x f_t+ E_t \cdot \nabla_v f_t = 0 \quad &\text{in}\quad(0,\infty)\times\T^d\times\R^d,\\
    \rho_t(x) = \int_{\R^d} f_t(x,v)\,\dv \quad &\text{in}\quad(0,\infty)\times\T^d,\\
     -\Delta_x \Phi[\rho_t] = \rho_t - e^{\Phi[\rho_t]}\quad &\text{in}\quad(0,\infty)\times\T^d,\\
    E_t = - \nabla_x \Phi[\rho_t] &\text{in}\quad(0,\infty)\times\T^d.\\
\end{cases}
\end{align}
In the above, we identify $\T^d = \lt[-\frac{1}{2},\frac{1}{2}\rt]^d$, and $f_t(x,v)$ represents the ion distribution in phase space $(x,v)\in \T^d \times\R^d$ at time $t$, while $\rho_t$ represents the spatial density of ions. The term $E_t=-\nabla_x \Phi_t[\rho_t]$ denotes the electrostatic force, where the electrostatic potential $\Phi_t[\rho_t]$ solves the Poisson--Boltzmann equation \eqref{iVP}$_{3}$. The exponential term $e^{\Phi_t}$ corresponds to the spatial density of massless electrons.

The IVP system formally arises by taking the electron mass to zero in the coupled two-species (electrons and ions) Vlasov--Poisson equations. Such a limit is justified by the physical observation that the mass of electrons is substantially lighter than that of ions.
At the timescale of the ions, electrons can be assumed to reach thermodynamic equilibrium (which corresponds to $e^{\Phi_t}$ in \eqref{iVP}$_{3}$), which is the reason the Poisson--Boltzmann equation (and thus the IVP system) arises. For a rigorous derivation of the IVP system from the two-component ion-electron model, we refer to \cite{bardosgolsenguyensentis2018,flynnguo2024} and the references therein. We also refer to \cite{ griffinpickeringiacobelli2021Torus, griffinpickeringiacobelli2021Recent,hankwan2011} for more
physical explanations on the IVP system.

%
%
%
%
%

Although the exponential nonlinearity makes the IVP system considerably more challenging than the classical Vlasov--Poisson (VP) system, significant progress has been made on the Cauchy problem for IVP in recent years. The existence of weak solutions to IVP in $\R^3$ was studied by F. Bouchut in \cite{bouchut1991}.  Using Wasserstein stability estimates, D. Han-Kwan and M. Iacobelli \cite{hankwaniacobelli2017a,hankwaniacobelli2017b} investigated the quasineutral limit for IVP in one and higher dimensions. Furthermore, M. Griffin-Pickering and M. Iacobelli studied the global well-posedness in $\T^d$ with $d=2,3$ \cite{griffinpickeringiacobelli2021Torus}, and in $\R^3$ \cite{griffinpickeringiacobelli2021Strong}. Later, L. Cesbron and M. Iacobelli studied the global well-posedness of IVP in bounded domains \cite{cesbroniacobelli2023}. See also \cite{griffinpickeringiacobelli2021Recent} for a broader survey of the well-posedness theory for \eqref{iVP}. 

In the context of global well-posedness theory, the initial data is typically assumed to have sufficient integrability with suitable velocity moments. For example, in \cite{griffinpickeringiacobelli2021Torus}, the initial data is assumed to satisfy $f_0 \in L^1 \cap L^\infty$ with specific velocity decay conditions, along with an $m$-th velocity moment bound for $m > 6$. As noted in \cite[Remark 1.5]{griffinpickeringiacobelli2021Strong}, these assumptions are essential, in the spirit of \cite{lionsperthame1991}, to ensure the boundedness of the local density $\rho$. This boundedness is in turn crucial when applying the approach of G. Loeper \cite{loeper2006}.
 
However, even when the data is not sufficiently integrable, it is reasonable to expect at least the existence of solutions, most typically in a distributional sense. In fact, even when the PDE does not make sense as distributions, it has been known since the seminal work of R. J. DiPerna and P. L. Lions that one may construct solutions which are physically meaningful. Indeed, for the classical VP system, DiPerna and Lions utilized the concept of renormalized solutions \cite{dipernalions1988,dipernalions1989a}, proving global existence under finite energy and initial data in the Orlicz space $L\log (1+L)(\R^d\times\R^d)$. The concept of renormalized solutions can be viewed as a generalization of distributional solutions from an Eulerian point of view.

Recently, in the context of transport equations, several results have emerged that connect the Eulerian and Lagrangian viewpoints. Indeed, L. Ambrosio et al. \cite{ambrosiocolombofigalli2017} proved the superposition principle for the classical VP system through a localized adaptation of the DiPerna--Lions theory \cite{dipernalions1989b}, and the equivalence of renormalized and Lagrangian solutions was established under very mild assumptions. In a similar spirit, but by a different argument building upon the results of \cite{ambrosiocrippa2008, BC13}, the existence of Lagrangian solutions to VP and their renormalized properties were obtained by A. Bohun et al. in \cite{BBC16b}. We refer also to \cite{BBC16,Xav18} for further developments.

Based on the above, we are naturally motivated to consider IVP in the low regularity setting, where $E_t\cdot \nabla_v f_t = \nabla_v\cdot(E_tf_t)$ in $\eqref{iVP}_1$ may lack meaning as a distribution. 
More concretely, let us consider the setup where we only know that $f_t \in L^1 \cap L^{q}$ for some $q>1$, and $|v|^2 f_t \in L^1$. Then an interpolation (see \eqref{eq: elem interpolation} below) shows that $\rho_t \in L^p$ for some $p>1$. Thus, the first natural question is whether the Poisson--Boltzmann equation $\eqref{iVP}_3$ is even solvable under only this integrability assumption on $\rho_t$. 

Toward this end, we recall the analysis in \cite{hankwaniacobelli2017a} where the Poisson--Boltzmann equation was tackled by decomposing it into a linear and nonlinear part. This same decomposition was also utilized in \cite{cesbroniacobelli2023,griffinpickeringiacobelli2021Torus}. Although only specific integrabilities and dimensions were considered in their works, following their arguments, one can in fact prove well-posedness of $\eqref{iVP}_3$ provided that the density $\rho_t$ lies in $L^p$ with $p>d/2$ for all $d\ge2$. We explain this in more detail in Section \ref{sec: d2}.


However, our goal is to handle the more general case of $p>1$ in all dimensions $d\ge2$, and it appears that the decomposition method of \cite{hankwaniacobelli2017a} is not so well-suited to this setting (see Remark \ref{rmk:cru} for more details). Therefore, the theory of well-posedness and stability for the Poisson--Boltzmann equation
\bq\label{eq:pb}
-\Delta \Phi[\rho] = \rho - e^{\Phi[\rho]} 
\eq
has to be significantly advanced by some other method. We introduce a novel decomposition which allows us to handle any $p>1$ and any dimension $d\ge 2$. In Section \ref{ssec:cov} this new decomposition is described in detail, and it is utilized to resolve the well-posedness of \eqref{eq:pb} based on the method of Calculus of Variations. Our method of proof uses crucially the damping effect of the exponential term $e^{\Phi[\rho]}$. Consequently, we are able to prove that for any $p>1$ and $d\ge 2$, the potential $\Phi[\rho]$ has $W^{1,2}(\T^d)$ regularity. 
The theorem reads as follows: 




\begin{theorem}[Well-posedness and properties of the Poisson--Boltzmann equation]
\label{thm:pb}
Let $ d \geq 2 $ and $ \rho \in L_+^p(\T^d) $ with $ p > 1 $, where
\[
L^p_+(\T^d):=\lt\{\rho\in L^p(\T^d): \rho(x) \ge 0 \,\,\textit{a.e.},\,\, \|\rho\|_{L^1(\T^d)}>0\rt\}.
\]
Then, the Poisson--Boltzmann equation \eqref{eq:pb} admits a unique solution $ \Phi[\rho] \in W^{1,2}(\T^d)$. It satisfies \bq\label{eq:Lr}
    \|e^{\Phi[\rho]}\|_{L^r(\T^d)}\le \|\rho\|_{L^r(\T^d)}\qquad\text{for all}\,\,\,\,r\in[1,p],
    \eq
    with equality when $r=1$ (commonly referred to as the neutrality condition), and following properties:
\begin{enumerate}[(i)]
    \item \textbf{$W^{1,2}$ estimates.} The solution $\Phi[\rho]$ satisfies
    \[
    \|\Phi[\rho]\|_{W^{1,2}}\lesssim B(p,\|\rho\|_{L^1}, \|\rho\|_{L^p}),
    \]
    where $B(p,a,b) \to \infty$ as either $p \to 1^+$, $a \to 0^+$, or $b \to +\infty$.
    \item \textbf{Weak Stability.} 
   For any $\rho \in L^p_+(\T^d)$, if $\rho_n \weakto \rho$ in $L^p(\T^d)$, then we have:
   \[
   \nabla \Phi[\rho_n] \weakto \nabla \Phi[\rho],\quad e^{\Phi[\rho_n]}\weakto e^{\Phi[\rho]}
   \] 
   in $L^2(\T^d)$ and $L^p(\T^d)$, respectively.
    
  \item \textbf{Uniform Lower Bound of $e^{\Phi[\rho]}$.} For any $\rho \in L^p_+(\T^d)$, we have
  \[
  e^{\Phi[\rho]} \ge c_\rho >0
  \]
  where $c_\rho$ depends on $d$, $\|\rho\|_{L^p}$, and $\|\rho\|_{L^1}$. As a result of \eqref{eq:Lr} and the lower bound estimate, we have
    \begin{equation} \label{eq: Lq}
       \Phi[\rho] \in L^q(\T^d)
 \end{equation}
for any $q \in [1,\infty)$.

     \item \textbf{Strong Stability.}  
    For any $\rho_1, \rho_2 \in L_+^p(\T^d)$ with $p >1$ and $q \in [1,\infty)$, the following stability estimate holds: 
    \bq\label{eq:str}
    \|\nabla \Phi[\rho_1] - \nabla \Phi[\rho_2]\|_{L^2}^2 + \lambda \|\Phi[\rho_1]-\Phi[\rho_2]\|_{L^q}^{ \frac{2q}{\min\{q,2\}} } \leq \Lambda \|\rho_1 - \rho_2\|_{L^p}^{\frac{2p}{\max\{p,2\}}},
    \eq
    where $ \lambda, \Lambda > 0 $ depend only on $ d, p, q, \|\rho_i\|_{L^1}, \|\rho_i\|_{L^p} $ for $ i = 1,2 $.
    
    \item \textbf{Continuity and semicontinuity of energy functionals.}  
    The potential energy and entropy functionals:
    \[
    \calP[\rho] := \int_{\T^d} \frac{1}{2} |\nabla \Phi[\rho]|^2 \, \dx, \quad \calS[\rho] := \int_{\T^d} \left(e^{\Phi[\rho]} (\Phi[\rho] - 1) + 1\right) \, \dx,
    \]
    are continuous under strong $L^p$-convergence of the density $\rho$, and lower semicontinuous under weak $L^p$-convergence.   
\end{enumerate}
\end{theorem}
\begin{remark}
    The estimates obtained in this work blow up as $p$ approaches $1$, and so we were unable to incorporate the case $p=1$. However, after the initial version of this work was completed, we learned that the work of Br\'ezis and Strauss \cite{BS73} proved well-posedness of general elliptic equations of the form $\calL \Phi + \beta(\Phi) = \rho$, in the limiting case where $\rho\in L^1$. Here $\calL$ is a second-order elliptic operator, $\beta$ a continuous monotone function, and the equation was considered in a bounded domain $\Omega$ with Dirichlet boundary conditions. We point out that in our setting of $\rho\in L^p_+$, $p>1$, adapting their result \cite[Corollary 12]{BS73} provides a unique solution in the class $\Phi\in W^{2,p}$. With only this regularity in hand, the integrability of $\Phi$ thus deteriorates as $p\to 1$. On the other hand, our main results obtain that $\Phi\in H^1\cap L^q$ for every $q\in [1,\infty)$, independently of the value of $p>1$. We are able to accomplish this by working directly on the energy functional associated to \eqref{eq:pb} in the framework of Calculus of Variations, and also by exploiting the exponential term on the right-hand side of \eqref{eq:pb}. Furthermore, the quantitative estimates in (iii)--(iv) appear to be genuinely new and unavailable in the previous literature. 
\end{remark}

\begin{remark}
    We remark that the space $L^p_+(\T^d)$ naturally aligns with the mathematical treatment of  \eqref{eq:pb}. While the non-negativity of $\rho$ has not been exploited in earlier works concerning \eqref{eq:pb}, it plays a critical role when studying the low integrability regime.
    
    Moreover, the requirement $\|\rho\|_{L^1(\T^d)}>0$ is not merely a technical condition but is essential for the existence of solutions. For instance, if $\rho\equiv m > 0$, the unique solution to \eqref{eq:pb} is explicitly written $\Phi = \log m$. As $m \to 0$, we get $\Phi \to -\infty$, demonstrating the need for the non-vanishing mass condition to exclude pathological cases.
\end{remark}

\begin{remark}[$W^{1,2}$ estimates and weak stability]

As aforementioned, observe that the constructed solution $\Phi[\rho]$ lies in $W^{1,2}(\T^d)$ for all $\rho \in L^p$ with $p>1$. This $L^2$-integrability of $\nabla \Phi[\rho]$ is stronger than what can be obtained typically via Calder\'on--Zygmund/Sobolev inequalities in the regime of small $p>1$. One may view this as a regularizing effect of the term $-e^\Phi$ in \eqref{eq:pb}.

In fact, in Corollary \ref{cor:w12b} we obtain quantitative $W^{1,2}$ estimates for $\Phi[\rho]$ with respect to  $\|\rho\|_{L^p}$ and  $\|\rho\|_{L^1}$. This yields, in particular, the weak stability results of Theorem \ref{thm:pb} (ii), which will be utilized to establish the energy inequality satisfied by solutions to IVP.
\end{remark}

\begin{remark}[Lower bound estimates and strong stability]
Strong $W^{1,2}$-type stability estimates \eqref{eq:str} are well-known at least for bounded densities \cite{griffinpickering2024,griffinpickeringiacobelli2021Torus}. These estimates take advantage of a positive lower bound of $e^{\Phi[\rho]}$. However, the methods used in \cite{griffinpickering2024} to derive the lower bound appear to require at least that $p>\frac{d}{2}$ (see Lemma \ref{lem:ePhi} where we recover a simpler proof). We then extend these results by proving a strict lower bound for $e^{\Phi[\rho]}$ whenever $\rho\in L^p_+$, $p>1$. Moreover a rather new approach is incorporated. Namely, we take advantage of the fact that solutions to the Poisson--Boltzmann equation enjoy a comparison principle. Our approach is motivated by the work of D. Lannes et al. \cite{lanneslinaressaut2013}, where the well-posedness of \eqref{eq:pb} was discussed for sufficiently regular densities which avoid vacuum. Note that our analysis does not require such assumptions.




\end{remark}

%
%
%
%
%

Now, we go back to our original motivation, which is the Cauchy problem for \eqref{iVP} under rough initial data. As aforementioned, we consider the following class of initial data:
\begin{align}
    \label{eq:inc}
    f_0 \ge 0, \quad f_0 \in L^1_2 \cap L^q(\T^d \times \R^d)\quad \text{for some }\,\,q>1, \quad \text{and} \quad \|f_0\|_{L^1(\T^d\times\R^d)} > 0,
    \end{align}
    where $L^1_2(\T^d\times\R^d):= L^1(\T^d\times\R^d,(1+|v|^2)\dx\dv)$. Assuming these bounds propagate in time (which will be possible in the finite energy class), the elementary interpolation inequality
\begin{align}
\label{eq: elem interpolation}
    \|\rho_t\|_{L^p(\T^d)} \lesssim \|f_t\|_{L^{q}(\T^d\times\R^d)}^{\frac{2/d}{1/q' + 2/d}} \||v|^2 f_t\|_{L^1(\T^d\times\R^d)}^{\frac{1/q'}{1/q' + 2/d}}, \qquad p := \frac{d(q-1)+2q}{d(q-1)+2}
\end{align}
gives $\rho_t \in L^\infty([0,\infty);L^p(\T^d))$ with $p>1$. Thus, incorporating the stability results of Theorem \ref{thm:pb}, we are able to construct global solutions to \eqref{iVP}:
    

\begin{theorem}\label{thm:lag}
Let $f_0$ satisfy \eqref{eq:inc}. Then there exists
\[
f_t \in C([0,+\infty);L^1(\T^d\times\R^d)) \cap L^\infty([0,\infty);L_2^1\cap L^q(\T^d\times\R^d))
\]
which is a Lagrangian (and renormalized) solution
to the ionic Vlasov--Poisson system \eqref{iVP} corresponding to initial data $f_0$ (see Definition \ref{def:lag}). Moreover, the following properties hold.
\begin{enumerate}[(i)]
\item \textbf{Energy inequality.}  For all $t \ge0$, the energy is bounded by its initial value:
\bq\label{eq:egyi}
\calE[f_t] \le \calE[f_0],
\eq
where
\bq\label{eq:egy}
\calE[f_t]:=\iint_{\T^d\times\R^d} \frac{1}{2}|v|^2 f_t \,\dx\dv + \intt \frac{1}{2}|\nabla_x\Phi[\rho_t]|^2\,\dx + \intt e^{\Phi[\rho_t]}(\Phi[\rho_t]-1)+1 \,\dx.
\eq
\item \textbf{Continuity of densities and electric field.} The ion density $\rho_t$ and the thermalized electron density $e^{\Phi[\rho_t]}$ are, as maps of time, strongly continuous in $L^r(\T^d)$ for any $r\in \left[1,\frac{d(q-1)+2q}{d(q-1)+2}\right)$, and the electric field $E_t$ is strongly continuous in $L^2(\T^d)$.
\end{enumerate}
\end{theorem}

The proof of this theorem is based largely on the framework developed in \cite{dipernalions1989a}. Approaching from the Eulerian point of view, the stability results of Theorem \ref{thm:pb} allow us to prove the existence of a renormalized solution to \eqref{iVP}. As stated in Theorem \ref{thm:lag}, this one is physical in the sense that it satisfies the energy inequality for every $t\ge 0$. Then we prove that this one is in fact a Lagrangian solution to \eqref{iVP} by means of the so-called extended superposition principle \cite[Theorem 4.10]{ambrosiocolombofigalli2017}. Although the results of \cite{ambrosiocolombofigalli2017} are stated mostly for the whole space $x\in \R^d$, the statements there are adaptable to the case of the flat torus; in order to avoid any ambiguity, we provide the necessary details and proofs in Appendix \ref{sec:pre}.

We close this section with some further remarks on Theorem \ref{thm:lag}.

\begin{remark}
 The initial energy $\calE[f_0]$ is finite for any initial data satisfying \eqref{eq:inc} (see Theorem \ref{thm:pb},(v)). 
\end{remark}

%
%
%
%
%

\begin{remark}
    By the classical theory of renormalized solutions, any renormalized solution to IVP is a weak solution as soon as its distributional formulation is well-defined. Of course, the latter holds as soon as $E_t f_t \in L^1_{\rm loc}(\R_+\times \T^d \times \R^d)$. Since $E_t \in L^\infty_tL^2_x$ (by Theorem \ref{thm:pb}), we see that $q \ge 2$ is a sufficient condition for global existence of distributional solutions, thereby extending the work of \cite{griffinpickeringiacobelli2021Torus} to arbitrary dimensions and lower integrability initial data. In fact, in the lower dimensions $d=2,3$ as considered in \cite{griffinpickeringiacobelli2021Torus}, the range of $q$ can be further relaxed: $q \ge \tfrac{1}{11}(12+3\sqrt{5}) \approx 1.701$ for $d=3$, and $q \ge \frac{1}{8}(7+\sqrt{17})\approx 1.390$ for $d=2$. This is because Sobolev embeddings provide higher integrability of $E_t$ in low dimensions, see \cite{dipernalions1988, dipernalions1988b} where the same thresholds appear for the classical VP system. 
  
\end{remark}

The rest of the paper is organized as follows. In Section \ref{sec:pb}, we study the well-posedness and stability of \eqref{eq:pb}. We then construct the global Lagrangian solutions to IVP in Section \ref{sec: lag}. The relevant results and notions of renormalized and Lagrangian solutions are recorded in Appendix \ref{sec:pre} for convenience of the reader.

%
%
%
%
%


%
%
%
%
%
\section{Analysis of the Poisson--Boltzmann Equation} \label{sec:pb}

In this section, we study the well-posedness and properties of the Poisson--Boltzmann equation for densities $\rho\in L^p_+(\T^d)$ with $p > 1$, thereby proving Theorem \ref{thm:pb}. We recall the equation:
\bq\label{eq:pb0}
-\Delta \Phi = \rho - e^{\Phi}
\eq
and we aim to establish solutions in $W^{1,2}(\T^d)$ with $d \ge2$. 

\subsection{Preliminaries: $p>\frac{d}{2}$ case}\label{sec: d2}

To set the stage, we briefly review prior results on \eqref{eq:pb0}, which hold under higher integrability assumptions on $\rho$. A key decomposition method for addressing was introduced in \cite{hankwaniacobelli2017b} for $d =1$ and extended in \cite{griffinpickeringiacobelli2021Torus} for $d \ge 2$. The decomposition involves splitting $\Phi$ as follows:
\bq\label{eq:pb_br}
-\Delta \overline{\Phi} = \rho - \intt \rho\, \dx
\eq
and
\bq\label{eq:pb_ht}    
-\Delta \widehat{\Phi} = \intt \rho \, \dx - e^{\overline{\Phi} + \widehat{\Phi}}.
\eq
Adding these equations, it is clear that 
\[
\Phi=\overline{\Phi}+\widehat{\Phi}
\]
solves \eqref{eq:pb0}. Using the Green's function $G$ for the negative Laplacian on the torus:
\[
-\Delta G := \delta_0 - 1\qquad \text{in}\,\,\,\calD'(\T^d),
\]
the solution to \eqref{eq:pb_br}, constrained to mean zero, is given explicitly by
\begin{align*}
    \overline{\Phi} = G * \rho.
\end{align*}

The following result provides regularity properties of $\overline{\Phi}$ under sufficiently high integrability of $\rho$, based on Calder\'on--Zygmund estimates and Sobolev embeddings \cite{gilbargtrudinger1977}. 
\begin{lemma}\label{lem:lin}
Let $d\ge 2$. For $\rho \in L^p_+(\T^d)$ with $\frac{d}{2}<p <d$, the solution $\overline{\Phi}$ to \eqref{eq:pb_br} satisfies
\[
\|\nabla \overline{\Phi}\|_{L^{p^*}(\T^d)} + \|\overline{\Phi}\|_{C^{0,\alpha}(\T^d)} \le N_{d,p}\lt\|\rho - \int_{\T^d}\rho\,\dx \rt\|_{L^p(\T^d)},
\]
where $p^*\in (d, \infty)$ is the Sobolev conjugate of $p$, and $\alpha= 2- \frac{d}{p} \in (0,1)$.
\end{lemma}

To solve \eqref{eq:pb0}, it remains to address the nonlinear equation \eqref{eq:pb_ht}. The following result was previously proven under specific regularity assumptions on $\overline{\Phi}$:
\begin{lemma}\cite[Proposition 3.5]{griffinpickeringiacobelli2021Torus}
\label{lem:IGP} 
Let $d\ge2$, and let $\overline{\Phi}$ satisfy
\begin{equation}\label{eq: assumption: igp}
\|\nabla\overline{\Phi}\|_{L^2(\T^d)} + \|\overline{\Phi}\|_{L^\infty(\T^d)} <+\infty.
\end{equation}
Then, for any given $m>0$, the equation
\begin{equation}\label{eq : m hat}
-\Delta \widehat{\Phi} = m - e^{\overline{\Phi}+ \widehat{\Phi}}
\end{equation}
has a unique solution $\widehat{\Phi} \in W^{1,2}(\T^d)$. If we further assume that $\overline{\Phi} \in C^{0,\alpha}(\T^d)$ for some $\alpha \in (0,1)$, then $\widehat{\Phi} \in C^{2,\alpha}(\T^d)$.
\end{lemma}

\begin{remark} In \cite{griffinpickeringiacobelli2021Torus}, the constant $m$ was typically chosen to be $1$. However, the result of Lemma \ref{lem:IGP} applies equally well for any $m>0$ since the proof is independent of this specific choice.
\end{remark}

\begin{remark}\label{rmk:cru}
The existence of $\widehat{\Phi}$  in Lemma \ref{lem:IGP} was obtained in \cite{griffinpickeringiacobelli2021Torus} via the Calculus of Variations. The conditions
\[
\overline{\Phi} \in L^\infty(\T^d)\quad\text{and}\quad m>0
\]
in \eqref{eq: assumption: igp} played a critical role in ensuring the compactness of minimizing sequences for the action functional associated to \eqref{eq : m hat}. In the case  of $\rho \in L^p(\T^d)$ with $p>\frac{d}{2}$, then $\overline{\Phi}$, the solution to \eqref{eq:pb_br}, satisfies the regularity condition \eqref{eq: assumption: igp}, as established in Lemma \ref{lem:lin}. Consequently, applying Lemma \ref{lem:IGP} with $m=\intt \rho\,\dx$ yields the unique solution $\widehat{\Phi}$ to \eqref{eq:pb_ht}, which yields the solution $\Phi=\overline{\Phi}+\widehat{\Phi}$ to \eqref{eq:pb0}. However, when $\rho \in L^p(\T^d)$ with $p \leq \frac d2$, $\overline{\Phi}$ is not necessarily bounded, which reveals a limitation of the decomposition method \eqref{eq:pb_br}--\eqref{eq:pb_ht} for densities with lower integrability. Addressing this challenge requires the development of new techniques.
\end{remark}

Throughout this section, we fix $p > 1$ and $d\ge2$ as the spatial dimension, unless stated otherwise.
%
%
%
%
%
\subsection{Solvability via Calculus of Variations: $p>1$ case}\label{ssec:cov}
To overcome the barriers of prior approaches based on the decomposition \eqref{eq:pb_br}--\eqref{eq:pb_ht}, we introduce a refined framework that decomposes both the density and the electric potential. This method accommodates the case where $\rho$ possesses low regularity, specifically $\rho \in L^p_+(\T^d)$ for $p > 1$.

For a given density $\rho\in L^p_+(\T^d)$, we define the following decomposition:
\begin{equation} \label{eq:rfl}
\begin{split}
    &\rho = \rho_{\flat} + \rho^\sharp, \quad \mbox{with} \quad \rho_\flat(x) :=  \rho(x) \wedge \lt(\frac{1}{2}\|\rho\|_{L^1(\T^d)}\rt),
    \end{split}
\end{equation}
where $\rho_{\flat}$ captures the bounded part of the density, and $\rho^\sharp = \rho - \rho_{\flat}$ represents the remainder. 

The corresponding electric potential is then decomposed into:
\[
\Phi = \Phi_\flat + \Phi^\sharp,
\]
with $\Phi_\flat $ and $\Phi^\sharp$ satisfying
\begin{align}
\label{eq:pb:fl}
    &-\Delta \Phi_\flat = \rho_\flat - \int_{\T^d} \rho_\flat \, \dx, \quad \int_{\T^d} \Phi_\flat = 0,\\
\label{eq:pb:sh}    
    &-\Delta \Phi^\sharp = \int_{\T^d} \rho_\flat \, \dx + \rho^\sharp - e^{\Phi_\flat + \Phi^\sharp}.
\end{align}
Since \eqref{eq:pb:fl} is linear and $\rho_\flat \in L^\infty(\T^d)$ by construction, the classical results guarantee the existence and uniqueness of $\Phi_\flat \in W^{1,2} \cap L^\infty(\T^d)$. Additionally, the mass term
\[
m_\flat:= \intt \rho_\flat\,\dx >0
\]
is strictly positive since $\rho$ is non-negative and not identically zero. As mentioned in Remark \ref{rmk:cru}, these properties allow us to reduce the nonlinear analysis to \eqref{eq:pb:sh}, where the primary challenge lies in handling the exponential term $e^{\Phi_\flat + \Phi^\sharp}$.

To address this challenge, we reformulate \eqref{eq:pb:sh} as a variational problem. The associated energy functional is:
\bq\label{calj}
h \mapsto \calJ[h]:= \intt \frac{1}{2}|\nabla h|^2 + \lt(e^{\Phi_\flat + h} - \lt(m_\flat + \rho^\sharp \rt)h\rt) \, \dx,\quad m_\flat:=\intt \rho_\flat\,\dx.
\eq
The Euler--Lagrange equation corresponding to $\calJ[\cdot]$ is precisely \eqref{eq:pb:sh}, ensuring that any minimizer is a solution to the nonlinear problem. Furthermore, strict convexity of $\calJ$ guarantees the uniqueness of the solution.

We now rigorously establish the existence and uniqueness of the solution $\Phi^\sharp$ to \eqref{eq:pb:sh} in the lemma below.

\begin{lemma}\label{lem:cov}
For any $\rho \in L_+^p(\T^d)$, there exists a unique solution $\Phi^\sharp \in W^{1,2}(\T^d)$ satisfying \eqref{eq:pb:sh} in the sense of distributions.
\end{lemma}
\begin{proof}
The proof is divided into four steps. First, we show that $\calJ[h]$ is bounded below and admits a minimizing sequence. Compactness arguments provide the convergence of this sequence to a weak limit, which we verify to be a minimizer of $\calJ[h]$. Finally, the minimizer is shown to satisfy \eqref{eq:pb:sh} in the sense of distributions.  

\medskip

\noindent \textbf{(Step I: Boundedness of the energy functional)} Consider the energy functional $\calJ$ appeared in \eqref{calj} in the class of periodic functions $h \in W^{1,2}(\T^d)$.  In this step, we aim to show here that the energy functional $\calJ$ is bounded below, and that therefore there exists a minimizing sequence of $\calJ$. By choosing $h=-\Phi_\flat \in W^{1,2}(\T^d) \cap L^\infty(\T^d)$, we first note that
\bq\label{eq:jub}
\begin{split}
\calJ[-\Phi_\flat] &= \intt  \frac{1}{2}|\nabla \Phi_\flat|^2 +1  +\lt(m_\flat + \rho^\sharp \rt)\Phi_\flat\,\dx\\
&= \intt \lt\{ \frac{1}{2}\lt( \rho_\flat - m_\flat\rt) + m_\flat + \rho^\sharp \rt\}\Phi_\flat +1 \,\dx \\
&\le 2\|\rho\|_{L^1}\|\Phi_\flat\|_{L^\infty} +1\\
&\le N_d\|\rho\|_{L^1}^2 +1 =:C_0.
\end{split}
\eq
For $h \in W^{1,2}(\T^d)$, using the decomposition $h= h_+- h_-$, where $h_+:=\max\{h,0\}$, we estimate lower bounds for the term in the functional $\calJ$. Observe that
\begin{equation}\label{eq:jfm}
\begin{split}
e^{\Phi_\flat + h} - \lt(m_\flat + \rho^\sharp \rt)h &=  e^{\Phi_\flat}e^{h}-\lt(m_\flat + \rho^\sharp \rt) h_+ + \lt(m_\flat + \rho^\sharp \rt) h_- \\
&\ge e^{-\|\Phi_\flat\|_{L^\infty}}e^{h}-m_\flat h_+ - \frac{(h_+)^{p'}}{p'}- \frac{\rho^p}{p} + \lt(m_\flat + \rho^\sharp \rt) h_- ,
\end{split}
\end{equation}
where $p' = \frac p{p-1}$ is the H\"older conjugate of $p$.

When $h<0$, it is clear that \eqref{eq:jfm} can be estimated further as
\begin{equation}
\label{eq: h<0}
\begin{split}
    e^{\Phi_\flat + h} - \lt(m_\flat + \rho^\sharp \rt)h & \ge e^{-\|\Phi_\flat\|_{L^\infty}}e^h + \lt(m_\flat + \rho^\sharp\rt) h_- - \frac{\rho^p}{p}  \ge \frac{1}{2}e^{-\|\Phi_\flat\|_{L^\infty}}e^h + m_\flat |h| + \rho^\sharp h_- - \frac{\rho^p}{p}.
    \end{split}
\end{equation}
On the other hand, when $h\ge 0$ we use the convexity of $x \mapsto e^x$ and its Taylor expansion to yield
\begin{equation}\label{eq:c1}
\begin{split}
e^{-\|\Phi_\flat\|_{L^\infty}}e^{h}-m_\flat h_+ - \frac{(h_+)^{p'}}{p'} &\ge  \frac{1}{2}e^{-\|\Phi_\flat\|_{L^\infty}}e^{h} + m_\flat h_+ -C_1 = \frac{1}{2}e^{-\|\Phi_\flat\|_{L^\infty}}e^h + m_\flat |h| - C_1,
\end{split}
\end{equation}
where 
\bq\label{eq:c1f}
C_1:= N_d(\|\rho\|_{L^1(\T^d)} + \|\rho\|_{L^1(\T^d)}^2 ) + (4e^{N_d\|\rho\|_{L^1(\T^d)}} \Gamma(p'+3))^{p'} >0
\eq
(see Lemma \ref{lem:C1} below for details). Consequently, for $h\ge 0$, \eqref{eq:jfm} can be estimated further as
\begin{align*}
    e^{\Phi_\flat+h} - \lt(m_\flat + \rho^\sharp \rt)h \ge \frac{1}{2}e^{-\|\Phi_\flat\|_{L^\infty}} e^h + m_\flat |h| - C_1 - \frac{\rho^p}{p}.
\end{align*}
Together with \eqref{eq: h<0}, we deduce that
\bq\label{eq:jlb}
\calJ[h] \ge \intt \frac{1}{2}|\nabla h|^2 + \frac{1}{2}e^{-\|\Phi_\flat\|_{L^\infty}}e^{h} +m_\flat |h| + \rho^\sharp h_-\,\dx -C_1 - \frac{1}{p}\intt \rho^p \,\dx
\eq
which confirms that $\calJ[\cdot]$ is bounded below.

\medskip

\noindent \textbf{(Step II: Compactness of minimizing sequence)} Let $h_k \in W^{1,2}(\T^d)$ be a minimizing sequence of $\calJ[\cdot]$:
\[
\lim_{k\to\infty} \calJ[h_k]=\inf_{h \in W^{1,2}(\T^d)} \calJ[h]=:\alpha \in \R.
\]
The boundedness of $\calJ$ established in {\bf (Step I)}, specifically, \eqref{eq:jub} and \eqref{eq:jlb}, implies that for all sufficiently large $k$, we have
\bq\label{eq:c2}
\intt \frac{1}{2}|\nabla h_k|^2 + m_\flat |h_k| \,\dx \le C_0+C_1 + \frac{1}{p}\intt \rho^p \,\dx  =:C_2< +\infty.
\eq
Note that we used the fact that $\rho^\sharp$ is non-negative. Since $m_\flat>0$, by Poincar\'e's inequality, $\{h_k\}$ is uniformly bounded in $L^2(\T^d)$ as well. Indeed:
\bq\label{eq:l2b}
\begin{split}
\|h_k\|_{L^2(\T^d)}^2 &= \lt\| h_k - \intt h_k \, \dx\rt\|_{L^2(\T^d)}^2 + \lt|\intt h_k \,\dx\rt|^2 \\
&\le N_d\|\nabla h_k\|_{L^2(\T^d)}^2 + \lt(\intt |h_k|\,\dx\rt)^2 \le 2N_dC_2 + \frac{C_2^2}{m_\flat^2}<+\infty.
\end{split}
\eq
These bound estimates yield that $h_k$ converges (up to a subsequence) weakly in $W^{1,2}(\T^d)$ to some limit $\Phi^\sharp \in W^{1,2}(\T^d)$:
\[
h_k \weakto \Phi^\sharp \qquad \text{in}\,\,\, W^{1,2}(\T^d) \quad \mbox{as} \quad k \to \infty.
\]
By the compact Sobolev embedding, up to a subsequence, we obtain
\[
h_k \to \Phi^\sharp \qquad \text{in}\,\,\, L^{2}(\T^d) \,\,\,\text{and}\,\,\,\text{a.e.} \quad \mbox{as} \quad k \to \infty.
\]
Further, by employing \eqref{eq:jub} and \eqref{eq:jlb}, we observe that
\[
\intt e^{h_k} \,\dx \le C_2 e^{\|\Phi_\flat\|_{L^\infty}}<+\infty,
\]
which implies that the $(h_k)_+$ are uniformly bounded in $L^r(\T^d)$ for each $r\in[1,+\infty)$. Consequently, by the almost everywhere convergence and the Vitali convergence theorem, we conclude that 
\begin{equation}
\label{eq: h_k+ conv}
(h_k)_+ \to \Phi^\sharp_+ \quad \text{in}\quad L^{p'}(\T^d) \quad \mbox{as} \quad k \to \infty.
\end{equation}

\medskip

\noindent \textbf{(Step III: $\Phi^\sharp$ is a minimizer)} We now confirm that the weak limit $\Phi^\sharp$ of the minimizing sequence $\{h_k\}$ is indeed a minimizer of the functional $\calJ$. 

By the lower semicontinuity of the $W^{1,2}$-norm under weak convergence:
\begin{equation}
\label{eq: conf1}
\liminf_{k \to \infty} \intt \frac{1}{2}|\nabla h_k|^2\, \dx \ge \intt \frac{1}{2}|\nabla \Phi^\sharp|^2\, \dx.
\end{equation}
Further, by Fatou's lemma, we get
\begin{equation}
\label{eq: conf2}
\liminf_{k \to \infty} \intt e^{\Phi_\flat + h_k}\, \dx \ge  \intt e^{\Phi_\flat+\Phi^\sharp}\, \dx.
\end{equation}
For the linear term, we separate it into contributions from the positive and negative parts of $h_k$. Since $\rho^\sharp \in L^p(\T^d)$, \eqref{eq: h_k+ conv} yields
\begin{equation}
\label{eq: conf: h+}
\lim_{k \to \infty} \intt -\lt(m_\flat + \rho^\sharp \rt) (h_k)_+ \, \dx = \intt -\lt(m_\flat + \rho^\sharp \rt) \Phi^\sharp_+ \, \dx.
\end{equation}
On the other hand, applying Fatou's lemma to $(h_k)_-$ gives
\begin{equation}
\label{eq: conf: h-}
\liminf_{k \to \infty} \intt  \lt(m_\flat + \rho^\sharp \rt) (h_k)_- \,\dx \ge  \intt  \lt(m_\flat + \rho^\sharp \rt) \Phi^\sharp_- \,\dx.
\end{equation}
Combining \eqref{eq: conf: h+} and \eqref{eq: conf: h-}, we deduce
\begin{equation}
\label{eq : conf3}
\begin{split}
    \liminf_{k\to\infty} \intt -\lt(m_\flat + \rho^\sharp \rt) h_k \,\dx &= \liminf_{k\to\infty} \left( -\intt \lt(m_\flat + \rho^\sharp \rt) (h_k)_+ \,\dx + \intt \lt(m_\flat + \rho^\sharp \rt) (h_k)_- \,\dx \right)\\
    &\ge \intt - \lt(m_\flat + \rho^\sharp \rt) \Phi^\sharp \,\dx.
    \end{split}
    \end{equation}
Substituting \eqref{eq: conf1}, \eqref{eq: conf2}, and \eqref{eq : conf3} into the expression for $\calJ$, we find
\bq\label{eq:jmm}
\alpha = \lim_{k\to \infty} \calJ[h_k] \ge \calJ[\Phi^\sharp],
\eq
proving that $\Phi^\sharp$ is a minimizer.

\medskip

\noindent\textbf{(Step IV: The minimizer satsifies the Euler--Lagrange equation)} We now verify that the minimizer $\Phi^\sharp$ satisfies the Euler--Lagrange equation \eqref{eq:pb:sh}. From the boundedness \eqref{eq:jlb} and the minimality of $\calJ[\Phi^\sharp]$ \eqref{eq:jmm}, we first note that
\[
\nabla \Phi^\sharp \in L^2(\T^d),\quad e^{\Phi^\sharp} \in L^{1}(\T^d),\quad \lt(1+\rho^\sharp\rt) \Phi^\sharp_- \in L^1(\T^d).
\]
For any $\phi \in C^\infty(\T^d)$ and $\eta>0$, the minimality of $\Phi^\sharp$ implies
\[
\calJ[\Phi^\sharp] \le \calJ[\Phi^\sharp + \eta \phi]. 
\]
This yields
\begin{align*}
0&\le \frac{\calJ[\Phi^\sharp + \eta \phi] - \calJ[\Phi^\sharp] }{\eta} \\
&= \frac{1}{\eta}\lt(\intt \frac{1}{2}|\nabla  \Phi^\sharp + \eta\nabla \phi|^2 - \frac{1}{2}|\nabla  \Phi^\sharp|^2\, \dx\rt) \\
&\quad + \frac{1}{\eta}\lt(\intt e^{\Phi_\flat}e^{ \Phi^\sharp+\eta \phi} -e^{\Phi_\flat}e^{ \Phi^\sharp}\,\dx \rt) + \frac{1}{\eta}\lt(\intt\, (m_\flat + \rho^\sharp) \lt(- \Phi^\sharp-\eta\phi +  \Phi^\sharp \rt)\dx\rt) \\
&= \intt \nabla \Phi^\sharp \cdot \nabla  \phi + \frac{\eta|\nabla \phi|^2}{2} \,\dx + \intt e^{\Phi_\flat+ \Phi^\sharp}\lt(\frac{e^{\eta\phi}-1}{\eta}\rt) \,\dx -\intt (m_\flat + \rho^\sharp)\phi\,\dx.
\end{align*}
By taking the limit as $\eta \to 0^+$, we find that
\[
0 \le \intt  \nabla\Phi^\sharp \cdot \nabla  \phi + e^{\Phi_\flat+ \Phi^\sharp} \phi - (m_\flat + \rho^\sharp)\phi\,\dx.
\]
Choosing $-\phi$ instead of $\phi$ gives us the opposite inequality, and therefore
\[
\intt  \nabla\Phi^\sharp \cdot \nabla  \phi + e^{\Phi_\flat+ \Phi^\sharp} \phi - (m_\flat + \rho^\sharp)\phi\,\dx=0.
\]
This establishes that $\Phi^\sharp$ satisfies \eqref{eq:pb:sh} in the sense of distributions.
\end{proof}

The constant $C_1$ in \eqref{eq:c1f} was introduced as part of the lower bound for the energy functional $\calJ$. While its definition guarantees that the functional is bounded below, it is beneficial to derive an explicit estimate for $C_1$ to provide further clarity and rigor in our analysis. The estimate is stated in the following lemma.

\begin{lemma}\label{lem:C1}
The constant $C_1$ in \eqref{eq:c1f} can be estimated as
\[
C_1 \le N_d(\|\rho\|_{L^1(\T^d)} + \|\rho\|_{L^1(\T^d)}^2 ) + (4e^{N_d\|\rho\|_{L^1(\T^d)}} \Gamma(p'+3))^{p'}.
\]
\end{lemma}
\begin{proof}
We aim to obtain the estimate for the constant $C_1>0$ in \eqref{eq:c1}, which satisfies
\begin{align*}
e^{-\|\Phi_\flat\|_{L^\infty}}e^{h}-m_\flat h - \frac{h^{p'}}{p'} &\ge  \frac{1}{2}e^{-\|\Phi_\flat\|_{L^\infty}}e^{h} + m_\flat h -C_1, \quad \forall \, h \geq 0.
\end{align*}

From the  convexity and positivity of $x \mapsto e^x$, we get
\[
e^h \ge e^s + e^s(h-s) \ge e^s(h-s) \quad \forall \, h,s\in \R.
\]
By choosing $s:=a +\log4 + \log(2b)$, where $a,b>0$, we find
\[
\frac{1}{4}e^{-a}e^h - bh \ge bh -(2b)(a+\log4 +\log(2b)).
\]
Now, choose $a:=\|\Phi_\flat\|_{L^\infty(\T^d)}\le N_d\|\rho\|_{L^1(\T^d)}$ and $b :=m_\flat \le \|\rho\|_{L^1(\T^d)}$. Then we obtain
\begin{equation}
\label{eq: poly}
\frac{1}{4}e^{-\|\Phi_\flat\|_{L^\infty}}e^h - m_\flat h \ge m_\flat h - N_d(\|\rho\|_{L^1(\T^d)} +\|\rho\|_{L^1(\T^d)}^2)
\end{equation}
for some constant $N_d>0$. On the other hand, we also note that
\begin{align}
\label{eq : h >>}
\frac{1}{4}e^{-a+h} - \frac{h^{p'}}{p'} \ge 0,\quad \forall \, h \ge 4e^a m!>1,\quad m:= \left\lceil p' \right\rceil +1.
\end{align}
Indeed, for such $h$, we deduce
\[
\frac{1}{4}e^{-a+h} = \frac{e^h}{4e^a} \ge \frac{h}{4e^am!}h^{m-1} \ge h^{p'} \ge \frac{h^{p'}}{p'}.
\]
In the opposite case where $0\le h< 4e^a m!$, we can simply estimate
\begin{equation*}
\begin{split}
\frac{1}{4}e^{-a+h} - \frac{h^{p'}}{p'}  \ge -\frac{h^{p'}}{p'}  \ge -(4e^a m!)^{p'} \ge - (4e^a\Gamma(p'+3))^{p'}.
\end{split}
\end{equation*}
This and \eqref{eq : h >>} imply
\begin{align*}
    \frac{1}{4}e^{-a+h} - \frac{h^{p'}}{p'} \ge -(4e^a \Gamma(p' + 3))^{p'} \quad \forall \, h\ge 0.
\end{align*}
Adding this with \eqref{eq: poly}, we obtain
\begin{align*}
    \frac{1}{2}e^{-a+h} - m_\flat h - \frac{h^{p'}}{p'} \ge m_\flat h - N_d(\|\rho\|_{L^1(\T^d)} + \|\rho\|_{L^1(\T^d)}^2 ) - (4e^{N_d\|\rho\|_{L^1(\T^d)}} \Gamma(p'+3))^{p'}.
\end{align*}
This completes the proof.
\end{proof}

As a consequence of Lemma \ref{lem:cov}, for each $\rho \in L^p_+(\T^d)$ we have the existence and uniqueness of a solution $\Phi\in W^{1,2}(\T^d)$ to \eqref{eq:pb0}. 

We now provide the bound estimate for the electron density as stated in \eqref{eq:Lr}.

\begin{lemma} Let $\Phi\in W^{1,2}(\T^d)$ be the unique solution  to \eqref{eq:pb0} constructed in Theorem \ref{thm:pb}. Then, we have
\bq
\label{eq:Lr2}
\begin{split}
&\|e^{\Phi[\rho]}\|_{L^1(\T^d)} = \|\rho\|_{L^1(\T^d)}, \quad \|e^{\Phi[\rho]}\|_{L^r(\T^d)}\le \|\rho\|_{L^r(\T^d)}\quad\text{for all}\,\,\,\,r\in(1,p].
\end{split}
\eq
\end{lemma}
\begin{proof}
For the case $r=1$, we use $\psi=1$ as a test function in \eqref{eq:pb0} to obtain
\begin{align*}
    0 = \int_{\T^d} (\rho - e^{\Phi})\,\dx.
\end{align*}
Since $\rho$ is integrable, the triangle inequality shows $e^{\Phi}$ is also integrable, and hence the result follows.

For $r\in (1,p]$, we consider $\psi = e^{(r-1)\Phi}$ as a test function in the weak formulation of \eqref{eq:pb0}. This leads to
\begin{align}\label{eq: ephi}
    \intt e^{r\Phi} \, \dx \le \intt \rho\, e^{(r-1)\Phi} \, \dx.
\end{align}
A formal calculation with H\"older's inequality yields the desired result immediately. To rigorously justify the calculation, we use the truncated test functions $\psi^k = e^{(r-1)(\Phi \wedge k)}$. Then we find
\begin{align} \label{eq: phik}
    \int_{\T^d} 1_{\{\Phi \le k\}} |\nabla \Phi|^2 e^{(r-1)(\Phi\wedge k)} \,\dx +\int_{\T^d} e^{\Phi} e^{(r-1)(\Phi\wedge k)}\,\dx = \int_{\T^d} \rho \, e^{(r-1)(\Phi \wedge k)}\,\dx, \quad \forall k\in \bbN. 
\end{align}
Since $\Phi\in H^1(\T^d)$ and $e^\Phi \in L^1(\T^d)$, all terms in this equation are well-defined. It follows by the density of smooth functions in $L^p(\T^d)$ spaces that \eqref{eq: phik} holds for every $k\in \bbN$, and in particular
\begin{align*}
    \int_{\T^d} e^{\Phi}e^{(r-1)(\Phi\wedge k)}\,\dx \le \int_{\T^d} \rho \, e^{(r-1)(\Phi\wedge k)}\,\dx, \quad \forall \, k\in \bbN.
\end{align*}
Taking the limit $k\to\infty$ and applying the monotone convergence theorem yields the desired bound estimates \eqref{eq: ephi}. 
\end{proof}

In fact, one can prove using \eqref{eq:Lr2} and Schauder estimates that under the additional assumption $\frac{d}{2}<p<d$, the solution $\widehat{\Phi}$ belongs to $C^{2,\alpha}(\T^d)$. This result was originally obtained in \cite[Proposition 3.1]{griffinpickeringiacobelli2021Torus}. Notably, our approach simplifies the proof and provides a sharper estimate, as the $C^{2,\alpha}$ norm of $\widehat{\Phi}$ is bounded by a constant that grows more moderately compared to the double-exponential dependence on $\|\rho\|_{L^p}$ established in \cite{griffinpickeringiacobelli2021Torus}.

\begin{lemma} Let $d\ge2$ and $\frac{d}{2}<p<d$. For any $\rho \in L^p_+(\T^d)$, the unique solution $\widehat{\Phi}$ to \eqref{eq:pb_br}--\eqref{eq:pb_ht} satisfies:
\begin{equation}\label{eq:htupplow}
    - N_{p,d} \|\rho\|_{L^p(\T^d)} + \log \|\rho\|_{L^1(\T^d)} \le \intt \widehat{\Phi} \,\dx \le \log \|\rho\|_{L^1(\T^d)},
\end{equation}
as well as
\[
\|\widehat{\Phi}\|_{C^{2,\alpha}(\T^d)} \le N_{p,d}\lt( \lt(1+ |\log \|\rho\|_{L^1(\T^d)}| + \|\rho\|_{L^p(\T^d)}\rt) \exp\lt({N_{p,d}\|\rho\|_{L^p(\T^d)}}\rt)\rt).
\]
Here $\alpha:= 2- \frac{d}{p} \in (0,1),$ and $N_{p,d}>0$ is a constant depending only on $p$ and $d$.
\end{lemma}
\begin{proof} Let $G$ denote the Green's function satisfying $-\Delta G = \delta_0 -1$ in $\calD'(\T^d)$. Then, the uniqueness part of Lemma \ref{lem:IGP} implies that the solution $\widehat{\Phi}$ can be represented as
\begin{equation}\label{eq3: widehatphi}
\widehat{\Phi} = \int_{\T^d} \widehat{\Phi} \,\dx + G * \lt(\intt \rho\,\dx - e^{\Phi}\rt).
\end{equation}

We begin by proving \eqref{eq:htupplow}, thus controlling the first term on the right-hand side of \eqref{eq3: widehatphi}. Since we already know from \cite{griffinpickeringiacobelli2021Torus} that $\widehat{\Phi}\in C^2$, we apply the maximum principle to \eqref{eq:pb_ht} and find that at a point of minimum $x_*\in \T^d$ of $\widehat{\Phi}$, it holds
\begin{align*}
    0 \ge \intt \rho \,\dx - e^{\overline{\Phi}(x_*) + \widehat{\Phi}(x_*)}.
\end{align*}
Consequently, we obtain the lower bound
\begin{align*}
    \widehat{\Phi} \ge - \|\ov{\Phi}\|_{L^\infty} + \log \lt(\intt \rho \,\dx \rt) \ge - N_{p,d}\|\rho\|_{L^p(\T^d)} + \log \lt(\intt \rho \,\dx \rt),
\end{align*}
the last inequality owing to Lemma \ref{lem:lin}.
For an upper bound, we proceed in the following way. Using that $\ov{\Phi}$ has mean zero, note that $\intt \widehat{\Phi} \,\dx = \intt \Phi\,\dx$. Furthermore, we know that $\intt \rho \,\dx = \intt e^{\Phi} \,\dx$. Therefore, Jensen's inequality shows
\begin{align*}
    \intt \widehat{\Phi} \,\dx = \intt \Phi \,\dx \le \log \lt(\intt \rho \,\dx \rt).
\end{align*}
Together these prove \eqref{eq:htupplow}.

We now prove the remaining assertion of the lemma. Since $\|e^\Phi\|_{L^p(\T^d)} \le \|\rho\|_{L^{p}(\T^d)}$, we find from \eqref{eq3: widehatphi} that:
\[
\|\nabla^2 \widehat{\Phi}\|_{L^p(\T^d)} \lesssim \|\rho\|_{L^1(\T^d)} + \|\rho\|_{L^p(\T^d)} \lesssim  \|\rho\|_{L^p(\T^d)},
\]
by Calder\'on--Zygmund estimates. Using $\|G\|_{L^1(\T^d)} + \|\nabla G\|_{L^1(\T^d)} <+\infty$ as well as \eqref{eq:htupplow}--\eqref{eq3: widehatphi}, we deduce
\[
\|\widehat{\Phi}\|_{W^{2,p}(\T^d)} \lesssim  |\log \|\rho\|_{L^1(\T^d)}|  + \|\rho\|_{L^p(\T^d)} . 
\]
Consequently, by Morrey's inequality, we find
\[
\|\widehat{\Phi}\|_{C^{0,\alpha}} \lesssim  |\log \|\rho\|_{L^1(\T^d)}|  + \|\rho\|_{L^p(\T^d)},
\]
where $\alpha= 2- \frac{d}{p}$. Next, we recall from Lemma \ref{lem:lin} that
\[
\|\overline{\Phi}\|_{C^{0,\alpha}} \lesssim  \|\rho\|_{L^p(\T^d)}.
\]
Since
\[
\|e^{g}\|_{C^{0,\alpha}(\T^d)} \le e^{\|g\|_{L^\infty(\T^d)}}\lt(1+\lt[g\rt]_{C^{0,\alpha}(\T^d)}\rt)
\]
for any $g \in C^{0,\alpha}(\T^d)$, we apply this to $\overline{\Phi} + \widehat{\Phi}$ to obtain
\[
\|e^\Phi\|_{C^{0,\alpha}(\T^d)} \lesssim \lt(1 + | \log \|\rho\|_{L^1} | +  \|\rho\|_{L^p(\T^d)}\rt)\exp\lt({N_{p,d}\|\rho\|_{L^p(\T^d)}}\rt).
\]
On the other hand, it is clear that the constant $\intt \rho\,\dx$ obeys
\[
\lt\|\intt \rho\,\dx\rt\|_{C^{0,\alpha}(\T^d)} = \|\rho\|_{L^1(\T^d)} \le  \|\rho\|_{L^p(\T^d)}.
\]
Hence, applying Schauder estimates to \eqref{eq3: widehatphi} concludes the desired result.
\end{proof}

Our work in this subsection completes the existence theory for the Poisson--Boltzmann equation \eqref{eq:pb0} when $\rho \in L^p_+(\T^d)$. In the sequel, we denote by
\bq\label{eq:Phi:def}
L^p_+(\T^d) \ni \rho \mapsto \Phi[\rho] \in W^{1,2}(\T^d)
\eq
the unique solution to \eqref{eq:pb0} in $W^{1,2}(\T^d)$.

%
%
%
%
%
\subsection{Weak stability}
In this subsection, we study the weak stability of the electrostatic potential $\Phi[\rho]$ with respect to the density $\rho$. The first step is to show that if a sequence of densities $\{\rho^n\}$ is bounded in $L^p_+(\T^d)$, then the corresponding sequence $\{\Phi[\rho^n]\}$ is bounded in $W^{1,2}(\T^d)$. Utilizing the inequalities \eqref{eq:c2} and \eqref{eq:l2b} allows us to yield this result (Corollary \ref{cor:w12b}), but we first need to investigate a lower bound for $m_\flat = \|\rho_\flat\|_{L^1(\T^d)}$.

\begin{lemma}\label{lem:mfl}
Let $d\ge 1$ and $p>1$. For $\rho \in L^p_+(\T^d)$ with $\rho_\flat$ as defined in \eqref{eq:rfl}, we have
\begin{equation}\label{eq:mfl}
m_\flat \ge \frac{1}{2}\|\rho\|_{L^1(\T^d)} \lt(\frac{\|\rho\|_{L^1(\T^d)}}{2\|\rho\|_{L^p(\T^d)}}\rt)^{\frac{p}{p-1}}>0.
\end{equation}
\end{lemma}
\begin{proof}
We establish the lower bound for $m_\flat$ by proving the following inequality:
\bq\label{eq:G}
\lt|\lt\{x \in \T^d : \rho(x) \ge \alpha \|\rho\|_{L^1(\T^d)}\rt\}\rt|  \ge \lt(\frac{(1-\alpha)\|\rho\|_{L^1(\T^d)}}{\|\rho\|_{L^p(\T^d)}}\rt)^{\frac{p}{p-1}},\quad \forall \,\alpha \in [0,1).
\eq
For $\alpha = 0$, \eqref{eq:G} follows directly from H\"older's inequality. For $\alpha \in (0,1)$, we define
\[
G_\alpha:= \lt\{x\in\T^d: \rho(x) \ge \alpha \|\rho\|_{L^1(\T^d)} \rt\},\quad B_\alpha:=\T^d\setminus G_\alpha.
\]
Using the definitions of $G_\alpha$ and $B_\alpha$, and applying H\"older's inequality, we obtain
\begin{align*}
\begin{aligned}
\|\rho\|_{L^1(\T^d)} 
= \left(\int_{B_\alpha} + \int_{G_\alpha}\right) \rho\,\dx \le  \alpha \|\rho\|_{L^1(\T^d)} + \|\rho\|_{L^p(\T^d)} |G_\alpha|^{1-\frac{1}{p}}.
\end{aligned}
\end{align*}
This gives
\[
|G_\alpha|^{1-\frac{1}{p}} \ge \frac{(1-\alpha)\|\rho\|_{L^1(\T^d)}}{\|\rho\|_{L^p(\T^d)}}.
\]
confirming \eqref{eq:G}. By definition of $\rho_\flat$, we have
\[
\rho_\flat(x) \ge \frac{1}{2}\|\rho\|_{L^1(\T^d)}1_{G_{\frac{1}{2}}}(x),
\]
and thus integrating both sides and applying \eqref{eq:G} with $\alpha=\frac{1}{2}$ yields the desired lower bound.
\end{proof}

This lower bound for $m_\flat$ shows that $\rho_\flat$ retains a positive contribution to the potential, even when $\rho$ has limited regularity. Using this result, we now derive a uniform bound for $\Phi[\rho]$ in $W^{1,2}(\T^d)$, provided the density satisfies certain integrability constraints.

\begin{corollary}\label{cor:w12b}
Let $d\ge 2$ with $p>1$. Suppose $\rho \in L^p_+(\T^d)$ satisfies the bounds: 
\begin{equation}
\label{eq: mM}
0<m \le \|\rho\|_{L^1(\T^d)} \le \|\rho\|_{L^p(\T^d)} \le M<+\infty.
\end{equation}
Then, we have
\begin{equation}
\label{eq : bdpmm}
\|\Phi[\rho]\|_{W^{1,2}(\T^d)} \le B(d,p,M,m)
\end{equation}
for some constant $B(d,p,M,m)>0$.
\end{corollary}
\begin{proof}
We estimate $\Phi[\rho]$ by separately considering its components $\Phi_\flat$ and $\Phi^\sharp$, which solve the linear and nonlinear subproblems, respectively. 

From the estimates in \eqref{eq:c2}--\eqref{eq:l2b}, we first find 
\[
\|\Phi^\sharp[\rho]\|_{W^{1,2}(\T^d)}^2 \le 4C_2 + \frac{C_2^2}{m_\flat^2}.
\]
The right-hand side has leading order (see Lemma \ref{lem:C1}, \eqref{eq:c2} and \eqref{eq:mfl})
\begin{align*}
    \left(\|\rho\|_{L^1(\T^d)} + \|\rho\|_{L^1(\T^d)}^2 + \Big(4e^{N_d \|\rho\|_{L^1(\T^d)}}\Gamma(p'+3)\Big)^{p'} +\frac{1}{p}\|\rho\|_{L^p(\T^d)}^p \right)^2 \cdot \frac{1}{\|\rho\|^2_{L^1(\T^d)}}\lt(\frac{\|\rho\|_{L^p(\T^d)}}{\|\rho\|_{L^1(\T^d)}}\rt)^{\frac{2p}{p-1}}.
\end{align*}
By inserting the bounds \eqref{eq: mM} into the above, we obtain \eqref{eq : bdpmm}. 

For the linear part, classical estimates for the Poisson equation yield
\[
\|\Phi_\flat[\rho]\|_{W^{1,2}(\T^d)}\le N_d\|\rho_{\flat}\|_{L^\infty(\T^d)} \le \frac{1}{2}N_d \|\rho\|_{L^1} \le \frac{1}{2}N_d M,
\]
where $N_d > 0$ is a constant depending on the dimension.

Combining these bounds for $\Phi_\flat$ and $\Phi^\sharp$ concludes the desired result.
\end{proof}

\begin{remark} Since the factor $\Gamma(p'+3)^{p'}$ grows unbounded as $p \downarrow 1$, equivalent to $p' \uparrow \infty$, the constant $B(d,p,M,m)$ appeared in Corollary \ref{cor:w12b} diverges to infinity as  $p \downarrow 1$.
\end{remark}

With this uniform bound on $\|\Phi[\rho]\|_{W^{1,2}(\T^d)}$, we are prepared to analyze the weak stability of the potential $\Phi[\rho]$ with respect to weak convergence of the density $\rho$. Specifically, the following lemma establishes the weak convergence of both $\nabla\Phi[\rho^n]$ and $e^{\Phi[\rho^n]}$ when the sequence $\{\rho^n\}$ converges weakly in $L^p(\T^d)$.

\begin{lemma}\label{lem:wst}
Suppose that $\rho_n \weakto \rho$ in $L^p(\T^d)$ for some $\rho \in L^p_+(\T^d)$. Then, we have
\begin{align}
\nabla\Phi[\rho^n] \weakto \nabla \Phi[\rho] \qquad &\mbox{in} \quad    L^2(\T^d), \label{eq:wc:ef}\\
e^{\Phi[\rho^n]} \weakto e^{\Phi[\rho]}  \qquad &\mbox{in} \quad  L^p (\T^d), \label{eq:wc:ed}
\end{align}
as $n \to \infty$.
\end{lemma}
\begin{proof}
We denote $\Phi_n:=\Phi[\rho^n]$ for brevity. Note that $\{\rho_n\}$ converges weakly also in $L^1(\T^d)$. Thus, by lower semicontinuity of the norm, $\|\rho_n\|_{L^1(\T^d)}$ is uniformly away from zero. Note also that $\|\rho_n\|_{L^p(\T^d)}$ is a bounded sequence. Hence, there exist constants $m,M>0$, independent of $n$, such that \eqref{eq: mM} holds for each $\rho_n$. Applying Corollary \ref{cor:w12b} shows $\{\Phi_n\}$ is bounded in $W^{1,2}(\T^d)$. Following the arguments in Lemma \ref{lem:cov}, there exists a subsequence $k_n$ and $\Phi \in W^{1,2}(\T^d)$ such that 
\begin{equation}
\label{eq:weakstab1}
\nabla \Phi_{k_n} \weakto \nabla \Phi\quad \mbox{in} \,\,\, L^2(\T^d)\quad \mbox{and} \quad \Phi_{k_n} \to \Phi \quad \text{a.e.} \quad \mbox{as} \quad n \to \infty.
\end{equation}
We claim that $\Phi = \Phi[\rho]$. Recall from Theorem \ref{thm:pb} that $\|e^{\Phi_n}\|_{L^p(\T^d)}\le \|\rho_n\|_{L^p(\T^d)}$. By Fatou's lemma, we confirm that $e^\Phi \in L^p(\T^d)$:
\[
\intt e^{p\Phi} \,\dx \le \liminf_{n \to \infty} \intt e^{p\Phi_{k_n}} \,\dx \le \liminf_{n \to \infty} \intt (\rho_{k_n})^p \,\dx \le C.
\]
This and \eqref{eq:weakstab1} imply
\[
e^{\Phi_{k_n}} \weakto e^{\Phi}  \quad\text{in}\quad L^p (\T^d) \quad \mbox{as} \quad n \to \infty.
\]
By passing to the limit $n \to \infty$, we observe that $\Phi$ verifies
\[
-\Delta \Phi = \rho - e^{\Phi} \quad\text{in}\quad \calD'(\T^d).
\]
Hence, the uniqueness of \eqref{eq:pb0} provided in Theorem \ref{thm:pb} shows $\Phi = \Phi[\rho]$. Since the limit is uniquely identified, our arguments above hold without having to pass to subsequences, and we have confirmed \eqref{eq:wc:ef}--\eqref{eq:wc:ed}, as desired.
\end{proof}

\begin{remark}\label{rmk:strcexp}
Thanks to the pointwise convergence in \eqref{eq:weakstab1} and the Vitali convergence theorem, there is in fact strong convergence of $e^{\Phi[\rho_n]}$:
\[
e^{\Phi[\rho^n]} \to e^{\Phi[\rho]} \quad\text{in}\quad L^r (\T^d),\quad \forall \, r \in [1,p) \quad \mbox{as} \quad n \to \infty.
\]
\end{remark}

%
%
%
%
%

\subsection{Comparison principle and uniform lower bound of $e^{\Phi[\rho]}$}

In this part, we establish a strict lower bound for $e^{\Phi[\rho]}$, which will be essential in deriving strong stability results. We begin by discussing a comparison principle, which asserts that $\Phi[\rho]$ inherits monotonicity from $\rho$ under pointwise ordering.

The comparison principle, stated informally, is as follows:
\[
\mbox{If}\,\,\, \rho_1(x) \le \rho_2(x) \,\,\,\mbox{for all}\,\,\, x \in \T^d,\,\,\,\mbox{then}\,\,\, \Phi[\rho_1](x) \le \Phi[\rho_2](x)\,\,\,\mbox{for all}\,\,\, x \in \T^d.
\]
The proof relies on a maximum/minimum principle argument, which we outline below. For this argument, $\Phi$ must possess sufficient regularity, specifically $C^2(\T^d)$. Thus, we temporarily assume that $\rho \in W^{1,r}(\T^d)$ with $r>d$. In this setting, it holds that $\Phi[\rho] \in C^{2,\alpha}(\T^d)$ for some $\alpha >0$, as follows from Morrey's inequality and Schauder estimates:
\[
\|\overline{\Phi}\|_{C^{2,\alpha}(\T^d)}\lesssim \|\rho\|_{C^{0,\alpha}(\T^d)} \lesssim \|\rho\|_{W^{1,r}(\T^d)}\quad \mbox{with } \alpha=1-\frac dr.
\]
Now, suppose $\rho_1, \rho_2 \in W^{1,r}(\T^d)$ (with $r>d$) satisfy
\[
 \rho_1(x) \le \rho_2(x) \,\,\,\mbox{for all}\,\,\, x \in \T^d.
\]
We observe that
\[
-\Delta (\Phi[\rho_2] - \Phi[\rho_1]) = \rho_2 - \rho_1 - \lt(e^{\Phi[\rho_2]} - e^{\Phi[\rho_1]}\rt).
\]
Let $x_\star$ be a point of minimum of $\Phi[\rho_2] - \Phi[\rho_1]$, i.e.,
\[
x_\star \in \argmin \lt(\Phi[\rho_2](\cdot) - \Phi[\rho_1](\cdot)\rt).
\]
Then, at $x_\star$, we get
\[
\Delta (\Phi[\rho_2] - \Phi[\rho_1])(x_\star) \ge 0,
\]
and this yields
\begin{align*}
e^{\Phi[\rho_2](x_\star)} - e^{\Phi[\rho_1](x_\star)} = \rho_2(x_\star) - \rho_1(x_\star) + \Delta (\Phi[\rho_2] - \Phi[\rho_1])(x_\star)\ge 0.
\end{align*}
Since the exponential function is monotonic, it follows that $\Phi[\rho_2](x_\star) - \Phi[\rho_1](x_\star) \ge 0$. By definition of $x_\star$, we deduce  
\[
\Phi[\rho_1](x) \le \Phi[\rho_2](x)\,\,\,\mbox{for all}\,\,\, x \in \T^d.
\]
Using this principle, we extend the result to densities in $L^p_+(\T^d)$ and derive a sharp uniform lower bound. To apply this principle to general $L^p_+(\T^d)$ densities, we need to remove the regularity assumption $\rho\in W^{1,r}(\T^d)$. The following lemma extends this result to $\rho_1, \rho_2 \in L^p_+(\T^d)$.
\begin{lemma}\label{lem:comp}
Suppose that $\rho_1, \rho_2 \in L^p_+(\T^d)$ satisfy
\[
\rho_1(x) \le \rho_2(x)\quad\textit{a.e.} \,\,\,x\in\T^d,
\]
then we have
\[
e^{\Phi[\rho_1](x)} \le e^{\Phi[\rho_2](x)}\quad\textit{a.e.}\,\,\, x\in\T^d.
\]
\end{lemma}
\begin{proof}
We consider a standard mollifier $\psi$ supported in $B_{1/4}(0)$. For each $n\in \N$, we set
\[
\psi_n(x):= n^d\psi(nx).
\]
We then define the mollified densities:
\[
\rho^n_i:= \rho_i * \psi_n,
\]
for $i=1,2$. Clearly, $\rho^n_i \in W^{1,r}(\T^d)$ for any $r>d$ and
\[
\rho_2^n - \rho_1^n = (\rho_2 - \rho_1)*\psi_n \ge 0.
\]
By the comparison principle, we have
\[
\Phi[\rho_1^n](x) \le \Phi[\rho_2^n](x)\,\,\,\mbox{for all}\,\,\, x \in \T^d.
\]
We then note that
\[
\|\rho^n_i\|_{L^1(\T^d)} = \|\rho_i\|_{L^1(\T^d)},\qquad \|\rho^n_i\|_{L^p(\T^d)} \le \|\rho_i\|_{L^p(\T^d)},
\]
and also
\[
\rho^n_i \to \rho_i \qquad \text{in}\,\,\, L^p(\T^d)\quad \text{as}\,\,\,n\to \infty.
\]
By Lemma \ref{lem:wst} and Remark \ref{rmk:strcexp}, we have 
\[
e^{\Phi[\rho_i^n]} \to e^{\Phi[\rho_i]}\qquad \text{in}\,\,\,L^1(\T^d).
\]
By passing to a subsequence, we obtain convergence a.e., and the lemma follows.
\end{proof}

Inspired by the arguments based on the maximum principle, we establish a method for obtaining a nontrivial lower bound of $e^{\Phi[\rho]}$. For the moment, we assume an additional condition of higher integrability on $\rho$, which we later remove to extend the result to the general $L^p_+$ setting. To streamline the proof, we borrow the decomposition method $\Phi = \overline{\Phi} + \widehat{\Phi}$ introduced in \cite{griffinpickeringiacobelli2021Torus} (see \eqref{eq:pb_br}--\eqref{eq:pb_ht}). Additionally, it is worth noting that the result presented here was obtained in 
the recent work \cite{griffinpickering2024} using a different approach. Our proof emphasizes simplicity and explicit estimates.

\begin{lemma}\label{lem:ePhi}
Let $d\ge2$ and $\rho\in L^p_+(\T^d)$ with $p >\frac{d}{2}$. Then $\Phi[\rho]$ satisfies:
\begin{equation}\label{eq:ePhi:lb}
e^{\Phi[\rho](x)} \ge \|\rho\|_{L^1(\T^d)}\exp\lt({-N_{d,p}\lt\|\rho -\intt \rho\,\dx\rt\|_{L^p(\T^d)}}\rt), \quad \forall x\in\T^d,
\end{equation}
where $N_{d,p}>0$ depends only on $d$, $p$.
\end{lemma}
\begin{proof}
Using the decomposition $\Phi = \overline{\Phi}(x) +\widehat{\Phi}(x)$, where $\widehat{\Phi} \in C^{2,\alpha}(\T^d)$, we can apply the minimum principle to the equation:
    \begin{align*}
        -\Delta \widehat{\Phi} = \intt \rho(x)\,\dx - e^{\overline{\Phi} + \widehat{\Phi}}.
    \end{align*}
By choosing $x_\star \in \argmin \widehat{\Phi} (\cdot)$, we obtain
\[
e^{\overline{\Phi}(x_\star) + \widehat{\Phi}(x_\star)} = \intt \rho(x)\,\dx+ \Delta \widehat{\Phi}(x_\star) \ge   \intt \rho(x)\,\dx.
\]
Thus, for any $x \in \T^d$, we have
\[
e^{\Phi(x)} = e^{\overline{\Phi}(x) +\widehat{\Phi}(x)} \ge e^{\overline{\Phi}(x) + \widehat{\Phi}(x_\star)}  = e^{\overline{\Phi}(x) -\overline{\Phi}(x_\star)}e^{\overline{\Phi}(x_\star) + \widehat{\Phi}(x_\star)} \ge e^{-2\|\overline{\Phi}\|_{L^
\infty(\T^d)}} \lt(\intt \rho \,\dx\rt).
\]
This combined with Lemma \ref{lem:lin} establishes the desired lower bound.
\end{proof}
\begin{remark}
Note that \eqref{eq:ePhi:lb} is sharp for constant densities; namely, we consider $\rho \equiv \bar \rho >0$, then we obtain explicit solution:
\[
\Phi[\bar \rho]=\log \bar \rho,
\]
so that \eqref{eq:ePhi:lb} becomes equality.
\end{remark}
\begin{remark}
The constant $N_{d,p}$ in Lemma \ref{lem:ePhi} satisfies
\[
N_{d,p} \to +\infty
\]
as $p \downarrow \frac{d}{2}$ so that the lower bound converges to $0$. This is due to the method of decomposition and the nature of Sobolev embeddings.
\end{remark}

To extend the lower bound to the setting $\rho \in L^{p}_+(\T^d)$ with $p>1$, we employ the comparison principle in Lemma \ref{lem:comp}. 

\begin{lemma}\label{lem:ulb}
Let $d\ge2$ with $p>1$ be given. For any $\rho\in L^p_+(\T^d)$, it holds that
\bq\label{eq:cush}
e^{\Phi[\rho](x)} \ge \frac{1}{2}\|\rho\|_{L^1(\T^d)} \lt(\frac{\|\rho\|_{L^1(\T^d)}}{2\|\rho\|_{L^p(\T^d)}}\rt)^{\frac{p}{p-1}} e^{-N_d\|\rho\|_{L^1(\T^d)}}\qquad \textit{a.e.}
\eq
for some $N_d>0$.
\end{lemma}

\begin{proof}
Since $\rho_\flat \le \rho$ \textit{a.e.} and $\rho_\flat \in L^\infty(\T^d)$, 
we apply Lemma \ref{lem:comp} and Lemma \ref{lem:ePhi} to observe that
\[
 e^{\Phi[\rho](x)} \ge   e^{\Phi[\rho_{\flat}](x)} \ge \|\rho_{\flat}\|_{L^1(\T^d)} e^{-N_d\|\rho_\flat - \intt \rho_{\flat} \, \dx\|_{L^\infty}} > 0\qquad \textit{a.e.}
\]
By the properties of $\rho_\flat$, we estimate
\[
\|\rho_\flat\|_{L^\infty} \le \frac{1}{2}\|\rho\|_{L^1(\T^d)},
\]
and also use the lower bound of $\|\rho_\flat\|_{L^1(\T^d)}$ in Lemma \ref{lem:mfl} to conclude the proof.
\end{proof}

\begin{remark}
In the case of the unit mass of density, i.e., $\|\rho\|_{L^1(\T^d)} =1$, we note that \eqref{eq:cush} provides a polynomial decay in $\|\rho\|_{L^p(\T^d)}$, whereas the bound in \eqref{eq:ePhi:lb} has exponential decay. In particular, for bounded $\rho \in L^\infty(\T^d)$ the previous lemma provides
\[
e^{\Phi[\rho](x)} \ge c_d\|\rho\|_{L^\infty(\T^d)}^{-1} \quad \text{for all}\,\,\, x\in\T^d.
\]
\end{remark}

As a corollary of this uniform lower bound for $e^{\Phi}$, we obtain the assertion \eqref{eq: Lq}:
\begin{corollary}\label{cor:Lq}
For each $\rho \in L^p_+(\T^d)$, it holds that
\[
\Phi[\rho] \in L^q(\T^d)
\]
for all $q \in [1,\infty)$.
\end{corollary}
\begin{proof}
We decompose $\Phi[\rho]$ into its positive and negative parts: $\Phi[\rho] = \Phi[\rho]_+ - \Phi[\rho]_-$. From Lemma \ref{lem:ulb}, we find
\[
\Phi[\rho]_- \in L^\infty(\T^d).
\]
Thus we only need check the integrability of the positive part. Since $e^{\Phi[\rho]} \in L^1(\T^d)$, Taylor expansion of the exponential shows that
$\Phi[\rho]_+ \in L^q(\T^d)$
for all $q\in [1,\infty)$.
\end{proof}

%
%
%
%
%
\subsection{Strong stability}
In this subsection, we present a stability estimate for the electric potential with respect to the $L^p$-norms of the densities.

\begin{lemma}
\label{lem:sst}
Let $d\ge2$ and $p>1$. For any $\rho_1, \rho_2\in L^p_+(\T^d)$ satisfying the bounds
\bq\label{eq:home}
m_1 \le \| \rho_i\|_{L^1} \le M_1, \quad \| \rho_i\|_{L^p} \le M_p,  \qquad i=1,2,
\eq
the following stability estimate holds:
\[
\intt |\nabla\Phi[\rho_1] - \nabla\Phi[\rho_2]|^2 \,\dx + \lambda \intt \lt(e^{|\Phi_1 -\Phi_2|}-1\rt)|\Phi_1 -\Phi_2| \,\dx \le \Lambda  \lt(\intt |\rho_1 - \rho_2|^p \,\dx\rt)^{\frac{2}{\max\{p,2\}}},  
\]
where $\lambda, \Lambda >0$ are constants depending only on $m_1,M_1,M_p$, and $d$. In particular, for any $q\ge 1$, it holds that
\begin{equation*}
\intt |\nabla\Phi[\rho_1] - \nabla\Phi[\rho_2]|^2 \,\dx + c_q{\lambda}  \lt(\intt |\Phi[\rho_1] - \Phi[\rho_2]|^{q} \,\dx\rt)^{\frac{2}{\min\{q,2\}}}  \le \Lambda \lt(\intt |\rho_1 - \rho_2|^p \,\dx\rt)^{\frac{2}{\max\{p,2\}}}  
\end{equation*}
for some $c_q>0$. 
\end{lemma}
\begin{proof}
We only need to discuss the case where $p\le 2$ and $q\ge2$. From this, the cases where $p>2$ or $q<2$ simply follow by an application of H\"older's inequality. We subtract the equations \eqref{eq:pb0} corresponding to $\rho_1$ and $\rho_2$. Corollary \ref{cor:Lq} authorizes the admission of $\Phi_1 - \Phi_2$ as a test function in the weak formulation, and we find
\[
\intt |\nabla \Phi_1 - \nabla \Phi_2|^2 \,\dx + \intt \lt(e^{\Phi_1}-e^{\Phi_2}\rt)(\Phi_1 -\Phi_2) \,\dx = \intt (\rho_1 - \rho_2)(\Phi_1 -\Phi_2)\,\dx.
\]
We note that
\[
\lt(e^{a}-e^{b}\rt)(a-b) \ge e^{\min\{a,b\}}(e^{|a-b|}-1)|a-b|,\quad \forall a,b \in\R.
\]
Due to Lemma \ref{lem:ulb}, both $e^{\Phi_1}$ and $e^{\Phi_2}$ have non-trivial uniform lower bounds. Namely, for some $c>0$ we have  
\bq\label{eq:ephc}
e^{\Phi_i(x)} \ge c>0 \quad \text{a.e.} \quad i=1,2.
\eq
Thus, it holds that
\begin{equation}
\label{eq: stab1}
\intt |\nabla \Phi_1 - \nabla \Phi_2|^2 \,\dx + c\intt \lt(e^{|\Phi_1 -\Phi_2|}-1\rt)|\Phi_1 -\Phi_2| \, \dx \le  \intt (\rho_1 - \rho_2)(\Phi_1 -\Phi_2)\,\dx.
\end{equation}
From the following inequality\footnote{By Taylor's expansion of the exponential function, we have
\[
\displaystyle (e^s-1)s =\sum_{k=2}^{\infty} \frac{s^k}{(k-1)!} \ge \frac{s^n}{\Gamma(n)}
\]
for any $2\le n\in \N$. When $q \in [2,\infty)\setminus \N$, deriving the sharp coefficient as in \eqref{eq:Gautschi} is not so trivial, but can be obtained through Gautschi's inequality.
}:
\bq\label{eq:Gautschi}
(e^{s}-1)s \ge \frac{s^q}{\Gamma(q)},\quad \forall s\ge0,\,\,\, \forall q \ge 2,
\eq
we get:
\begin{equation}
\label{eq: stab2}
\frac{c}{2}\intt \lt(e^{|\Phi_1 -\Phi_2|}-1\rt)|\Phi_1 -\Phi_2| \, \dx \ge \frac{c}{2\Gamma(q)} \intt |\Phi_1 -\Phi_2|^q\,\dx \quad \forall q\ge 2.
\end{equation}
On the other hand, we apply Young's inequality to estimate the right-hand side of \eqref{eq: stab1} as
\begin{equation}
\label{eq: stab3}
 \intt (\rho_1 - \rho_2)(\Phi_1 -\Phi_2) \,\dx \le \intt \frac{A^p}{p}|\rho_1-\rho_2|^p + \frac{1}{p'A^{p'}}|\Phi_1-\Phi_2|^{p'}\, \dx,
\end{equation}
which holds for any $A>0$. In view of \eqref{eq: stab2} (which holds with $q=p'$), if we choose $A>0$ to satisfy
\[
\frac{1}{p'A^{p'}} =\frac{c}{2\Gamma(p')},
\]
then the final term in \eqref{eq: stab3} can be absorbed to the left-hand side of \eqref{eq: stab1}. We deduce
\[
\intt |\nabla \Phi_1 - \nabla \Phi_2|^2 \,\dx + \frac{c}{2}\intt \lt(e^{|\Phi_1 -\Phi_2|}-1\rt)|\Phi_1 -\Phi_2|  \, \dx \le \frac{1}{p}\lt(\frac{2\Gamma(p')}{cp'}\rt)^{p-1}\intt |\rho_1 - \rho_2|^p \,\dx.
\]
By invoking the inequality \eqref{eq:cush}, we note that
\[
c=\frac{m_1}{2}\lt(\frac{m_1}{2M_p}\rt)^{\frac{p}{p-1}}e^{-N_d M_1} >0
\]
satisfies \eqref{eq:ephc}. Lastly, we take 
\[
\lambda :=\frac{c}{2},\quad \Lambda := \frac{1}{p}\lt(\frac{2\Gamma(p')}{cp'}\rt)^{p-1}
\]
to conclude the proof.
\end{proof}

%
%
%
%
%
\subsection{Energy functionals}

In this final subsection, we apply the previously obtained results to study the energy functionals introduced in \eqref{eq:egy}. 

We begin by recalling the electric potential energy and entropy of the thermalized electrons:
\[
\calP[\rho]= \intt \frac{1}{2}|\nabla\Phi[\rho]|^2 \,\dx \quad \mbox{and} \quad  \calS[\rho]:= \intt e^{\Phi[\rho]}(\Phi[\rho]-1)+1 \,\dx.
\]

The first result shows that these functionals are continuous under the strong topology of $L^p_+(\T^d)$.

\begin{lemma}
\label{lem:ect:str}
    Assume that $\{\rho_n\}\subset L^p_+(\T^d)$ is a sequence such that
\begin{align*}
\rho_n \to \rho \qquad \text{in}\quad L^p_+(\T^d) \quad \mbox{as} \quad n \to \infty.
\end{align*}
Then, we have
\[
\displaystyle \calP[\rho] = \lim_{n \to \infty} \calP[\rho_n],\qquad \calS[\rho]=\lim_{n \to \infty} \calS[\rho_n].
\]
\end{lemma}
\begin{proof}
We denote $\Phi_n:=\Phi[\rho_n]$ and $\Phi:=\Phi[\rho]$ for brevity. Due to the assumption of strong $L^p$ convergence of $\rho_n$, there exist $m_1,M_1,M_p>0$ (independent of $n$) such that \eqref{eq:home} is satisfied for each $\rho_n$ (and for $\rho$). Note that the convergence $\calP[\rho_n]\to \calP[\rho]$ thus follows from Lemma \ref{lem:sst}. On the other hand, Corollary \ref{cor:Lq} permits the admission of $\Phi_n$ as a test function in the weak formulation of \eqref{eq:pb0} (with $\rho_n$ and $\Phi_n$), which yields
\[
\intt |\nabla \Phi_n|^2\,\dx = \intt \rho_n \Phi_n\,\dx - \intt e^{\Phi_n}\Phi_n\,\dx.
\]
Rearranging the terms, we can rewrite
\begin{align*}
\calS[\rho_n] &=  \intt \rho_n \Phi_n\,\dx  - 2\calP[\rho_n] - \intt e^{\Phi_n}\,\dx+1 \\
&= \intt \rho_n \Phi_n\,\dx - 2\calP[\rho_n] - \intt \rho_n\,\dx + 1,
\end{align*}
the last line following by invoking the neutrality condition \eqref{eq:Lr}. Hence,
\begin{equation} \label{eq: sub S}
\calS[\rho_n]  - \calS[\rho] =  \intt (\rho_n \Phi_n -\rho \Phi )\,\dx - 2(\calP[\rho_n] -\calP[\rho]) -\intt (\rho_n -\rho) \,\dx.
\end{equation}
Applying Lemma \ref{lem:sst} and recalling Corollary \ref{cor:Lq}, we have that
\[
\|\rho_n \Phi_n -\rho \Phi \|_{L^1(\T^d)} \le \|\rho_n\|_{L^p(\T^d)}\|\Phi_n-\Phi\|_{L^{p'}(\T^d)} + \|\rho_n-\rho\|_{L^p(\T^d)} \|\Phi\|_{L^{p'}(\T^d)} \to 0,
\]
as $n$ tends to infinity. The convergence of the other two terms on the right-hand side of \eqref{eq: sub S} is immediate, hence $\calS[\rho_n] \to \calS[\rho]$, completing the proof.
\end{proof}

Next, we establish lower semicontinuity of the functionals $\calP$ and $\calS$ under weak topology of $L^p_+(\T^d)$.

\begin{lemma}\label{lem:lsc}
Assume that $\rho_n \weakto \rho$ in $L^p(\T^d)$ for some $\rho \in L^p_+(\T^d)$. Then we have
\begin{equation*}
\calP[\rho] \le \liminf_{n \to \infty} \calP[\rho_n] \qquad \mbox{and}\qquad \calS[\rho] \le \liminf_{n \to \infty} \calS[\rho_n].
\end{equation*}
\end{lemma}
\begin{proof}
The first assertion follows directly from Lemma \ref{lem:wst} and the lower semicontinuity of the $L^2$-norm under weak convergence.

For the entropy $\calS[\cdot]$, recall the convexity of the function $x \mapsto x(\log x -1) +1$. For any $a,b \geq 0$ (with $0 \log 0 =0$):
\[
a(\log a-1) + 1 \ge  b(\log b -1) +1 +\log b(a-b). 
\]
By taking $a=e^{\Phi_n(x)}$, $b=e^{\Phi(x)}$ and integrating over $x\in \T^d$, we obtain
\[
\calS[\rho_n] = \intt e^{\Phi_n}(\Phi_n-1)+1 \,\dx \ge \intt e^{\Phi}(\Phi-1)+1 \,\dx + \intt \Phi(e^{\Phi_n} - e^{\Phi})\,\dx =:\calS[\rho] + \calR_n
\]
From Corollary \ref{cor:Lq}, we know that $\Phi \in L^{p'}(\T^d)$. Employing Lemma \ref{lem:wst} (namely \eqref{eq:wc:ed}) shows $\calR_n \to 0$ as $n \to \infty$, and this completes the proof.
\end{proof}

%
%
%
%
%

%
%
%
%
%

\section{Construction of Solutions} \label{sec: lag}
In this section, we employ the results of Section \ref{sec:pb} to prove Theorem \ref{thm:lag} confirming the existence of Lagrangian solutions to \eqref{iVP}.

\subsection{Regularized system and a priori estimates}
We begin with an approximate model of \eqref{iVP}, introduced in \cite{griffinpickeringiacobelli2021Torus}, where the electric field is doubly-regularized. Consider a scaled mollifier 
\[
\chi_n(x):= n^d\chi(nx),\quad n \in \N,
\]
where $\chi:\T^d \to \R$ is a smooth, radially symmetric, non-negative function with unit mass, supported in a ball of radius $\frac{1}{4}$. The regularized model is given by
\begin{align}
\label{iVP:reg}
\begin{cases}
    \p_t f_t^n  + v \cdot \nabla_x f_t^n+ E_t^n \cdot \nabla_v f_t^n = 0 \quad &\text{in}\quad(0,\infty)\times\T^d\times\R^d,\\
    \rho_t^n(x) = \intr f_t^n(x,v)\,\dv \quad &\text{in}\quad(0,\infty)\times\T^d,\\
     -\Delta_x \Phi_t^n = \chi_n * \rho_t^n - e^{\Phi_t^n}\quad &\text{in}\quad(0,\infty)\times\T^d,\\
    E_t^n = - \chi_n * \nabla_x \Phi_t^n &\text{in}\quad(0,\infty)\times\T^d, \\  
\end{cases}
\end{align}
with the initial condition
\bq\label{iVP:reg:ini}
f_0^n|_{t=0} = f^n_0,
\eq
where $f_0^n \in C_c^\infty(\T^d\times\R^d)$ satisfies
\bq\label{eq:f0n}
f_0^n \to f_0\qquad \text{in}\quad L^1\cap L^q(\T^d\times\R^d)\quad \mbox{and}\quad \iint_{\T^d\times\R^d} |v|^2f_0^n \,\dx\dv \to \iint_{\T^d\times\R^d} |v|^2f_0 \,\dx\dv.
\eq
We remark that $\eqref{iVP:reg}_1$ can be written in the form of a continuity equation
\begin{equation*}
\begin{split}
    \p_t f_t^n + \textnormal{div}_{x,v}(\bfb_t^n f_t^n) = 0 \quad \mbox{with} \quad  \bfb_t^n(x,v) := (v,E_t^n(x)).
    \end{split}
\end{equation*}

The existence of solutions to the regularized system \eqref{iVP:reg}--\eqref{iVP:reg:ini} is well-established in \cite{griffinpickeringiacobelli2021Torus}, based on a fixed point argument in the Wasserstein distance. Due to the regularization, $E_t^n$ and $\nabla_x E_t^n$ are bounded in $[0,\infty)\times \T^d$, and $\bfb_t^n=(v,E_t^n)$ is a Lipschitz vector field on $\T^d\times \R^d$. As a result, the flow $\bfX^n(t):\T^d\times\R^d\to\R^{2d}$ associated to $\bfb_t^n$ is well-defined and unique. By classical results for the transport equation, we have 
\begin{equation*}
f_t^n = f_0^n \circ \bfX^n(t)^{-1}.
\end{equation*}
Since $\bfb_t^n$ is divergence-free, the flow is incompressible, ensuring
\[
\|\rho^n_t\|_{L^1(\T^d)} = \|f_t^n\|_{L^1(\T^d\times\R^d)} = \|f_0^n\|_{L^1(\T^d\times\R^d)},
\]
and 
\bq\label{eq:lqc:reg}
 \|f_t^n\|_{L^q(\T^d\times\R^d)} = \|f_0^n\|_{L^q(\T^d\times\R^d)}.
\eq

A key advantage of the double-regularization approach is that, as illustrated in \cite{horst1990}, the following energy functional is conserved in time:
\bq\label{ceg:reg}
\calE^n[f_t^n]:=\iint_{\T^d\times\R^d} \frac{1}{2}|v|^2f_t^n\,\dx\dv + \intt \frac{1}{2} |\nabla_x \Phi_t^n|^2\,\dx + \intt \lt[ e^{\Phi_t^n}(\Phi_t^n-1)+1 \rt]\,\dx.
\eq
We point out that this provides nice uniform bounds for $f^n$ because
\bq\label{ciegy}
\calE^n[f_0^n] \to \calE[f_0]\qquad \mbox{as}\quad n \to \infty.
\eq
Indeed, in the definition of $\calE^n[f_0^n]$ as written in \eqref{ceg:reg}, the convergence of the first term follows directly from \eqref{eq:f0n}. To handle the remaining terms, note that from the perspective of \eqref{eq:Phi:def}:
\[
\Phi_t^n = \Phi[\chi_n*\rho_t^n], \qquad E_t^n = - \chi_n * \nabla_x \Phi[\chi_n*\rho_t^n].
\]
By an elementary interpolation (see \eqref{eq: elem interpolation}), setting $p:= \frac{d(q-1)+2q}{d(q-1) +2}>1$, we observe that \eqref{eq:f0n} implies
\[
\rho^n_0 \to \rho_0\qquad \mbox{in}\quad L^p(\T^d) \quad \mbox{as}\quad n \to \infty.
\]
Since
\begin{align*}
\|\chi_n * \rho_0^n - \rho_0\|_{L^p(\T^d)} &\le \|\chi_n * (\rho_0^n - \rho_0)\|_{L^p(\T^d)} + \|\chi_n * \rho_0 - \rho_0\|_{L^p(\T^d)}  \cr
&\le \|\rho_0^n - \rho_0\|_{L^p(\T^d)} + \|\chi_n * \rho_0 - \rho_0\|_{L^p(\T^d)},
\end{align*}
we obtain
\begin{equation}
\label{eq: conviniteng1}
\chi_n * \rho_0^n \to \rho_0\qquad \mbox{in}\quad L^p(\T^d) \quad \mbox{as}\quad n \to \infty.
\end{equation}
Noting that $\|\rho_0\|_{L^1(\T^d)} = \|f_0\|_{L^1(\T^d\times\R^d)} >0$, it follows from Lemma \ref{lem:ect:str} that the last term in \eqref{ceg:reg} converges:
\[
 \intt e^{\Phi_0^n}(\Phi_0^n-1)+1\,\dx \to \intt \lt[ e^{\Phi[\rho_0]}(\Phi[\rho_0]-1)+1 \rt] \,\dx\qquad \mbox{as} \quad n \to \infty.
\]
On the other hand, from \eqref{eq: conviniteng1} and Lemma \ref{lem:sst}, we have
\[
\nabla_x \Phi_0^n = \nabla_x \Phi[\chi_n*\rho_0^n] \to \nabla_x \Phi[\rho_0]\qquad \mbox{in}\quad L^2(\T^d)
\]
as $n\to \infty$. Consequently, it is clear that the second term in \eqref{ceg:reg} converges:
\[
\intt \frac{1}{2} |\nabla_x \Phi_0^n|^2\,\dx \to  \intt \frac{1}{2} |\nabla_x \Phi[\rho_0]|^2 \,\dx
\]
as $n\to \infty$. This proves the convergence of initial energy \eqref{ciegy} as desired.

Next, we point out the $L^\infty([0,\infty);L^1(\T^d\times\R^d))$ norm of $|v|^2f$. By the conservation of energy and \eqref{ciegy}, we have
\[
\iint_{\T^d\times\R^d} \frac{1}{2}|v|^2f_t^n\,\dx\dv \le \calE^n[f_t^n]=\calE^n[f_0^n] \le C
\]
since the other two terms in \eqref{ceg:reg} are nonnegative. Hence, 
\bq\label{eq:v2m:uni}
\sup_{n\in \N}\sup_{t\in[0,\infty)}\iint_{\T^d\times\R^d}|v|^2f_t^n \,\dx\dv \le C.
\eq

%
%
%
%
%

\subsection{Passing to the limit in the electric field}

In this section we shall incorporate the results of Theorem \ref{thm:pb} in order to pass to the limit in \eqref{iVP:reg}. Thanks to \eqref{eq:lqc:reg}--\eqref{eq:v2m:uni} and the Dunford--Pettis criterion there is a limit $f\in L^\infty_{\rm loc}([0,\infty); L^1\cap L^q(\T^d\times \R^d))$ such that $f^n \weakto f$ weakly in $L^1((0,T)\times \T^d \times \R^d)$ for every $T>0$.

\begin{lemma}\label{lem:rho:strongconv}
    We set 
    \[
\rho_t:=\intr f_t \,\dv.
\]
Then for every $T>0$, it holds up to some subsequence that
\begin{equation*}
    \rho^n \to \rho \quad \text{strongly in} \quad L^{\ov{p}}((0,T)\times \T^d)
\end{equation*}
for all $\ov{p}\in [1,p)$ with $p := \frac{d(q-1)+2q}{d(q-1)+2}$.
\end{lemma}
\begin{proof}
We consider for each $\delta>0$, as in \cite{dipernalions1989a}, the sequence $\beta_\delta(f^n) := f^n/(1+\delta f^n)$. Since $f^n$ is a solution to \eqref{iVP:reg} in the classical sense of differential calculus, $\beta_\delta(f^n)$ solves
\begin{equation}\label{eq:betadelta}
    \p_t \beta_\delta(f^n) + v\cdot \nabla_x \beta_\delta(f^n) + E_t^n \cdot \nabla_v \beta_\delta(f^n) = 0.
\end{equation}
By applying velocity averaging lemmas \cite{DLM91}, we know that the sequence of averages $\lt\{\int_{\R^d} \beta_\delta(f^n) \varphi(v) \,\dv \rt\}_n$ is compact in $L^1((0,T)\times\T^d)$ for each $\delta>0$ and $\varphi\in \calD(\R^d)$. To show that such a compactness property holds for $f^n$ (instead of $\beta_\delta(f^n)$), let us assume without loss of generality that $q\in (1,2)$ in \eqref{eq:lqc:reg}; then we have
\begin{equation} \label{eq: beta app}
\begin{split}
    \int_0^T \iinttr |\beta_\delta(f^n) - f^n| \,\dx\dv\dt      &= \int_0^T \iinttr \frac{ \delta (f^n)^2 }{1 + \delta f^n} \,\dx\dv\dt \\
    &= \int_0^T \iinttr \frac{\delta^{q-1} (f^n)^q  }{(1+\delta f^n)^{q-1}} \lt(\frac{ \delta f^n }{1+\delta f^n}\rt)^{2-q} \,\dx\dv\dt \\
    &\le \int_0^T \iinttr \frac{\delta^{q-1} (f^n)^q  }{(1+\delta f^n)^{q-1}} \,\dx\dv\dt \\
    &\le \delta^{q-1} \int_0^T \iinttr (f^n)^q \,\dx\dv\dt \le C_T \delta^{q-1}.
\end{split}
\end{equation}
This estimate shows that for each $\varphi\in \calD(\R^d)$ the sequence $\lt\{\int_{\R^d} f^n \varphi(v) \,\dv \rt\}_n$ is also strongly compact in $L^1((0,T)\times\T^d)$. At last, using \eqref{eq:v2m:uni}, we deduce that this compactness property holds for every $\varphi\in C^0(\R^d)$ satisfying $|\varphi(v)|/|v|^2 \xrightarrow[|v|\to\infty]{}0$. In particular, choosing $\varphi(v)\equiv 1$, we obtain that
\begin{equation*}
    \rho^n \to \rho \quad \text{strongly in} \quad L^1((0,T)\times\T^d).
\end{equation*}
Finally, by interpolating the $L^1_2 \cap L^q$ bounds for $f^n$ as in \eqref{eq: elem interpolation}, we know that $\{\rho^n\}$ is bounded in $L^\infty(\R_+;L^p(\T^d))$. Therefore, by $L^p$-interpolation, we conclude
\begin{equation*}
    \rho^n \to \rho \quad \text{strongly in} \quad L^{\ov{p}}((0,T)\times \T^d), \quad \forall \ov{p} \in \lt[1, p \rt), \quad p := \frac{d(q-1)+2q}{d(q-1)+2}.
\end{equation*}
\end{proof}

As an application of Lemma \ref{lem:rho:strongconv} and Theorem \ref{thm:pb} we may prove the

\begin{lemma}\label{lem:elecstr}
    Let $\Phi[\rho_t]$ denote the unique solution to the Poisson--Boltzmann equation corresponding to $\rho_t$ (see Lemma \ref{lem:rho:strongconv} and Theorem \ref{thm:pb}). Then up to some further subsequence, we have
    \begin{equation*}
        E_t^n \to E_t := - \nabla_x \Phi[\rho_t] \quad \text{a.e. and in} \quad L^2_{\rm loc}(\R_+\times \T^d).
    \end{equation*}
\end{lemma}

\begin{proof}
We first apply Fubini's theorem to deduce strong convergence of the temporal traces $\rho_t^n$ as follows. From Lemma \ref{lem:rho:strongconv}, we know that for every $T>0$,
\begin{equation*}
        \int_0^T \lt(\intt |\rho_t^n - \rho_t|^{\ov{p}} \,\dx \rt) \dt \xrightarrow[n\to\infty]{} 0.
    \end{equation*}
	Thus, passing to some subsequence, it holds for almost every $t\in [0,T]$ that
    \begin{equation*}
        \intt |\rho_t^n - \rho_t|^{\ov{p}} \,\dx \xrightarrow[n\to\infty]{} 0.
    \end{equation*}
Taking $T\to\infty$ and extracting diagonally, we get a subsequence for which $\rho_t^n \to \rho_t$ in $L^{\ov{p}}(\T^d)$ for a.e. $t\ge 0$. Consequently, for a.e. $t\ge 0$:
\begin{equation*}
    \chi_n * \rho_t^n - \rho_t = \chi_n * (\rho_t^n - \rho_t) + (\chi_n  * \rho_t - \rho_t) \to 0 \quad \text{in} \quad L^{\ov{p}}(\T^d), \quad \forall \ov{p}\in [1,p).
\end{equation*}
As a result, applying Theorem \ref{thm:pb}, we obtain that for a.e. $t\ge 0$,
\begin{equation*}
	\|\nabla_x \Phi[\chi_n * \rho_t^n] - \nabla_x \Phi[\rho_t]\|_{L^2(\T^d)} = \|\nabla_x \Phi_t^n - \nabla_x \Phi[\rho_t]\|_{L^2(\T^d)} \to 0.
\end{equation*}
We can then easily check that this implies
\begin{equation*}
    \|E_t^n - E_t\|_{L^2(\T^d)} \to 0 \quad \text{a.e. }t\ge 0.
\end{equation*}
On the other hand, since $\{\chi_n * \rho_t^n\}$ is bounded in $L^\infty(\R_+;L^p(\T^d))$, Theorem \ref{thm:pb} tells us that $\{E_t^n\}$ is bounded in $L^\infty(\R_+;L^2(\T^d))$. In particular, $\lt\{\|E_t^n - E_t \|_{L^2(\T^d)}\rt\}_n$ is uniformly bounded in $L^\infty(\R_+)$. Therefore, the dominated convergence theorem of Vitali shows
\begin{align*}
	\forall T>0, \quad \int_0^T \lt\|E_t^n - E_t \rt\|_{L^2(\T^d)}^2 \dt \to 0.
\end{align*}
The thesis of the lemma follows by diagonal extractions.
\end{proof}

%
%
%
%
%

\subsection{Proof of Theorem \ref{thm:lag}}

\subsubsection{Existence of a renormalized/Lagrangian solution}

We follow largely the scheme of \cite{dipernalions1989a}. Taking again $\beta_\delta(f^n) := f^n/(1+\delta f^n)$, we recall that this one is a solution to \eqref{eq:betadelta} in the classical sense.

Now we shall pass to the limit $n\to\infty$ in the equation \eqref{eq:betadelta}. Letting $f_\delta$ denote the weak$^*$ limit of $\{\beta_\delta(f^n)\}_n$ in $L^\infty(\R_+\times \T^d \times \R^d)$, with the aid of Lemma \ref{lem:elecstr} we can easily pass to the limit to obtain that $f_\delta$ solves
\begin{equation}\label{eq:fdel}
    \p_t f_\delta + v\cdot \nabla_x f_\delta + E_t\cdot \nabla_v f_\delta = 0 \quad \text{in} \quad \calD'([0,\infty)\times \T^d\times\R^d),
\end{equation}
with initial condition $f_{\delta 0} = \beta_\delta(f_0)$. Then thanks to the facts that $f_\delta\in L^\infty(\R_+\times\T^d\times\R^d)$ and $E_t\in L^\infty(\R_+;L^2(\T^d))$, the results of DiPerna--Lions \cite{dipernalions1989b} apply, showing that $f_\delta$ satisfies at least for every $\beta\in C^1\cap W^{1,\infty}(\R_+)$ that
\begin{equation}\label{eq:betafdel}
    \p_t \beta(f_\delta) + v\cdot \nabla_x \beta(f_\delta) + E_t \cdot \nabla_v \beta(f_\delta) = 0 \quad \text{in} \quad \calD'([0,\infty)\times \T^d\times\R^d).
\end{equation}

We then claim that $f_\delta \to f$ strongly in a suitable space; since $\beta\in W^{1,\infty}(\R_+)$, this would allow us to pass to the limit $\delta\to 0$ in the above, showing that $f$ is a renormalized solution to \eqref{iVP}. Toward this end, we note that the estimate in \eqref{eq: beta app} actually provides (see also \eqref{eq:lqc:reg})
\begin{equation}\label{eq:carbonara}
    \sup_{t\in [0,T]} \|\beta_\delta(f^n) - f^n \|_{L^1(\T^d\times \R^d)} \le C\delta^{q-1},
\end{equation}
where the constant is independent of $n$. This yields:

\begin{lemma}\label{lem:carbonara}
    For each $T>0$, there holds
    \begin{equation*}
        f_\delta \in C([0,T];L^1(\T^d\times\R^d))
    \end{equation*}
    as well as
    \begin{equation*}
        f_{\delta} \to f \quad \text{in} \quad C([0,T];L^1(\T^d\times\R^d)).
    \end{equation*}
\end{lemma}
\begin{proof}
    First, we remark that we can easily check $f_\delta\in C([0,T];L^1(\T^d\times\R^d))$ from the equation \eqref{eq:fdel} that it satisfies. For instance, let us mollify both sides of \eqref{eq:fdel} by using a mollifier in the $(x,v)$-variables. Then the mollified version $f_\delta^\ve$ lies in $C([0,T];L^1_{\rm loc}(\T^d\times\R^d))$ and satisfies
    \begin{equation*}
        \p_t f_\delta^\ve + v\cdot \nabla_x f_\delta^\ve + E_t\cdot \nabla_v f_\delta^\ve = r^\ve
    \end{equation*}
    in the classical sense of differential calculus, where $\ve>0$ is the regularization parameter, and $r^\ve$ is the commutator term satisfying at least $r^\ve\to 0$ in $L^1_{\rm loc}(\R_+\times \T^d\times\R^d)$ (see \cite{dipernalions1989b}). Then, we consider the equation satisfied by $(f_\delta^\ve - f_\delta)$. Owing to \cite{dipernalions1989b}, $(f_\delta^\ve - f_\delta)$ is again a renormalized solution to this equation. We may therefore consider the equation satisfied by $\beta(f_\delta^\ve - f_\delta)$, for $\beta\in C^1\cap W^{1,\infty}(\R)$:
    \begin{equation*}
        \p_t \beta(f_\delta^\ve - f_\delta) + v\cdot \nabla_x \beta(f_\delta^\ve - f_\delta) + E_t \cdot \nabla_v \beta(f_\delta^\ve - f_\delta) = r^\ve \beta'(f_\delta^\ve - f_\delta).
    \end{equation*}
    Let us assume specifically that $\beta(0)=0$. Then, writing out the Green's formula that is satisfied by $\beta(f_{\delta t}^\ve - f_{\delta t})$, and using the fact that $f_\delta^\ve\to f_\delta$ strongly in $L^p_{\rm loc}(\R_+\times\T^d\times\R^d)$ for every $p\in [1,\infty)$, we may prove that for any test function $\varphi(x,v)\in C_c^\infty(\T^d\times\R^d)$, it holds that
    \begin{equation*}
        \forall t\in [0,T], \quad \iinttr \beta(f_{\delta t}^\ve - f_{\delta t}) \varphi(x,v) \,\dx\dv = o(1) \quad \text{as} \quad \ve\to 0,
    \end{equation*}
    with bounds on the right-hand side depending on $\|\beta'\|_{L^\infty}$ but not on $t\in [0,T]$.
    Then, we take a sequence $\{\beta_\kappa(z)\}_{\kappa\in \bbN}$ with $\beta_\kappa(z) \nearrow |z|$ and $\|\beta_\kappa'\|_{L^\infty} = O(1)$. From the estimate above, we can pass to the limit $\kappa\to\infty$ to obtain that
    \begin{equation*}
        \sup_{t\in [0,T]} \iinttr |f_{\delta t}^\ve - f_{\delta t}| \varphi(x,v) \,\dx\dv = o(1) \quad \text{as} \quad \ve \to 0.
    \end{equation*}
    In other words $f_{\delta}^\ve \to f_{\delta}$ in $C([0,T];L^1_{\rm loc}(\T^d\times\R^d))$ and so we have obtained that $f_\delta \in C([0,T];L^1_{\rm loc}(\T^d\times\R^d))$. Since $f_\delta$ lies in $L^\infty(0,T;L^1_2(\T^d\times\R^d))$ (remembering $\beta_\delta(f^n)\le f^n$ and \eqref{eq:v2m:uni}), we conclude to $f_\delta\in C([0,T];L^1(\T^d\times\R^d))$.

    Next, we consider arbitrary test functions $0\le \psi(t)\in L^1\cap L^\infty(0,T)$ and $\eta(x,v)\in L^\infty(\T^d\times\R^d)$. Let us further assume that $\|\psi\|_{L^1(0,T)}\le 1$ and $\|\eta\|_{L^\infty(\T^d\times\R^d)}\le 1$. Then thanks to \eqref{eq:carbonara} and $\|\eta\|_{L^\infty}\le 1$,
    \begin{align*}
        &\int_0^T \iinttr (\beta_\delta(f^n) - f^n) \psi(t) \eta(x,v) \,\dx\dv\dt \le C \delta^{q-1} \int_0^T \psi(t)\,\dt.
    \end{align*}
    Since $\psi\eta \in L^\infty((0,T)\times\T^d\times\R^d)$, we can pass to the limit $n\to\infty$ and obtain
    \begin{equation*}
        \int_0^T \iinttr (f_\delta - f) \psi(t) \eta(x,v) \,\dx\dv\dt \le C\delta^{q-1} \int_0^T \psi(t) \,\dt. 
    \end{equation*}
    Taking supremum over all such $\psi$ and $\eta$, we deduce
    \begin{equation*}
        \|f_\delta - f\|_{L^\infty(0,T;L^1(\T^d\times\R^d))} \le C\delta^{q-1} \xrightarrow[\delta\to 0]{} 0.
    \end{equation*}
    In the above, we used the fact that $L^1\cap L^\infty(0,T)$ is dense in $L^1(0,T)$. Since $f_\delta \in C([0,T];L^1)$, the above implies
    \begin{equation*}
        \|f_\delta - f\|_{C([0,T];L^1(\T^d\times\R^d))} \xrightarrow[\delta\to 0]{} 0,
    \end{equation*}
    as desired.
\end{proof}

Thanks to Lemma \ref{lem:carbonara}, we can now straightforwardly pass to the limit in \eqref{eq:betafdel} to obtain that $\beta(f)$ satisfies \eqref{iVP} in the sense of $\calD'([0,\infty)\times \T^d\times\R^d)$ for every $\beta\in C^1\cap W^{1,\infty}(\R_+)$, and next for every $\beta\in C^1\cap L^\infty(\R_+)$, by an approximation argument. In other words $f$ is a renormalized solution to \eqref{iVP}.

Finally, from Lemma \ref{lem:carbonara} we also know that $f\in C([0,\infty);L^1(\T^d\times\R^d))$, and therefore we may apply the superposition principle in Proposition \ref{prop: re are lag}. Therefore, we obtain that $f$ is a Lagrangian solution to \eqref{iVP}. This completes the proof of the existence of a global solution to \eqref{iVP}.


\subsubsection{Continuity of the maps $f_t$, $\rho_t$, $e^{\Phi[\rho_t]}$, and $E_t$}
We now complete the proof of assertion (ii) of Theorem \ref{thm:lag}. Specifically, we verify the temporal continuity of the maps $f_t$, $\rho_t$, $e^{\Phi[\rho_t]}$, and $E_t$ as follows:
\begin{enumerate}[(i)]
\item We know that the solution satisfies $f_t \in C([0,\infty);L^1(\T^d \times \R^d))$ and has bounded second moments. Using standard interpolation techniques, we find that the map $t\mapsto \rho_t$ is continuous in $L^{r}(\T^d)$ for all $1\le r< p= \frac{d(q-1)+2q}{d(q-1)+2}$.
\item In view of Lemma \ref{lem:wst}, more precisely, the remark following it, we conclude that $t\mapsto e^{\Phi[\rho_t]}$ is strongly continuous in $L^r(\T^d)$ for all $1\le r < p$.
\item By employing the strong stability provided in Lemma \ref{lem:sst}, we immediately deduce that the map $t\mapsto E_t$ is strongly continuous in $L^2(\T^d)$.
\end{enumerate}

%
%
%
%
%

\subsubsection{Energy inequality for a.e. $t\ge 0$}
We now verify the energy inequality \eqref{eq:egyi}, completing the proof of assertion (i) of Theorem \ref{thm:lag}. Let $\phi \in C_c((0,\infty))$ be a nonnegative, bounded test function. Since $f^n \weakto f$ in $L^1((0,T)\times\T^d\times\R^d)$, we can easily check that
\[
\int_0^\infty \phi(t)\iint_{\T^d\times\R^d} \frac{1}{2}|v|^2f_t\,\dx\dv\dt \le \liminf_{n \to \infty} \int_0^\infty \phi(t)\iint_{\T^d\times\R^d} \frac{1}{2}|v|^2f_t^n\,\dx\dv\dt,
\]
for instance by using appropriate cutoffs in the $v$-variable and then passing to the limit.

Next, by the proof of Lemma \ref{lem:elecstr}, we have
\[
\int_0^\infty \phi(t) \intt
\frac{1}{2}|\nabla_x\Phi[\rho_t]|^2\,\dx\dt
=
\lim_{n\to \infty}
\int_0^\infty \phi(t)\intt
\frac{1}{2}|\nabla_x\Phi_t^n|^2\,\dx\dt,
\]
thanks to the fact that $\phi$ is compactly supported. Finally, by the convexity of $x \mapsto x(\log x -1) +1$ and also its non-negativity, we have
\[
\phi(t)\lt(e^{\Phi_t^n}(\Phi_t^n - 1) + 1\rt) \ge  \phi(t)\lt(e^{\Phi[\rho_t]}(\Phi[\rho_t]-1) + 1\rt) +\phi(t) \Phi[\rho_t]\lt(e^{\Phi_t^n}-e^{\Phi[\rho_t]}\rt). 
\]
We focus on the last term here. On the one hand, by repeating the proof of Lemma \ref{lem:elecstr}, we know that $\Phi^n\to \Phi[\rho_t]$ a.e. and in $L^2_{\rm loc}(\R_+\times \T^d)$. On the other hand, by Theorem \ref{thm:pb} (more specifically Corollary \ref{cor:Lq}), we have $\Phi[\rho_t]\in L^\infty(\R_+;L^r(\T^d))$ for every $r\in [1,\infty)$. In particular $\phi(t) \Phi[\rho_t] \in L^r(\R_+\times \T^d)$ for every $r\in [1,\infty)$. Therefore, that last term vanishes in $L^1(\R_+\times \T^d)$ as $n\to\infty$, and we deduce
\[
\int_0^\infty \phi(t)\intt \lt[e^{\Phi[\rho_t]}(\Phi[\rho_t]-1) + 1\rt]\,\dx\dt  \le \liminf_{n\to \infty} \int_0^\infty \phi(t)\intt \lt[e^{\Phi_t^n}(\Phi_t^n-1) + 1\rt]\,\dx\dt. 
\]
Adding these inequalities together and using \eqref{ceg:reg}--\eqref{ciegy}, we deduce
\[
\int_0^\infty \phi(t) \calE[f_t]\,\dt \le \lt(\int_0^\infty \phi(t) \,\dt\rt) \calE[f_0].
\]
Since $\phi(t)$ is arbitrary, we deduce that
\begin{equation}
\label{eq: egyi ae}
\calE[f_t] \le \calE[f_0]\qquad \mbox{for a.e.}\,\,\, t \in [0,\infty).
\end{equation}
%
%
%
%
%
\subsubsection{Extension of the energy inequality to all $t\ge0$}
The kinetic energy (second moments) is lower semicontinuous under convergence in $L^1_{loc}(\T^d\times\R^d)$. Additionally, from Lemma \ref{lem:lsc}, the electric potential energy $\calP[\rho]$ and the entropy of the thermalized electrons $\calS[\rho]$ are lower semicontinuous with respect to the weak topology of $L^p(\T^d)$ for any $p>1$.

From the results of the previous section, the map $t\mapsto f_t$ is strongly continuous in $L^1(\T^d\times\R^d)$, and $t\mapsto \rho_t$ is strongly continuous in $L^r(\T^d)$ for all $1\le r < \frac{d(q-1)+2q}{d(q-1)+2}$. Note that Lemma \ref{lem:lsc} can be applied with any appropriate $r \in \left(1,\frac{d(q-1)+2q}{d(q-1)+2}\right)$ instead of $p$.

For any fixed $t_*\in (0,\infty)$, choose a sequence $t_n \to t_*$ such that \eqref{eq: egyi ae} holds for $t=t_n$. Then, by the lower semicontinuity properties discussed above, we deduce that \eqref{eq: egyi ae} holds with $t=t_*$. Since $t_*$ is arbitrary, this argument establishes that \eqref{eq: egyi ae} holds for all $t>0$, completing the proof of the theorem.

%
%
%
%
%
%
%
%
%

\section*{Acknowledgments}
This work is supported by NRF grants no. 2022R1A2C1002820 and RS-2024-00406821. DK is supported by NRF grant no. RS-2025-02312778.

\bigskip
\noindent\textbf{Conflict of Interest}

The authors declare that there is no conflict of interest regarding the publication of this paper.

\bigskip
\noindent\textbf{Data Availability}

No data were used to support this study.

\appendix
\section{Review of Renormalized and Lagrangian solutions} \label{sec:pre}

In this appendix, we review some of the notions and results which are used throughout this work. First, we state precisely our notion of a renormalized solution. Next, we recall the theory developed in \cite{ambrosiocolombofigalli2015, ambrosiocolombofigalli2017, ambrosiocrippa2008} and more specifically the superposition principle, which allows us to assert that the renormalized solution we constructed is indeed a Lagrangian solution to \eqref{iVP}. We provide the adaptations required in order to prove the superposition principle in our setting, as the formulations in \cite{ambrosiocolombofigalli2017} were originally stated only for $x\in \R^d$ (whereas we are in the periodic domain). Although the adaptations are relatively straightforward, we give all relevant terminology and results/proofs where deemed appropriate, for the convenience of the reader. For a more comprehensive discussion of these notions, we refer to \cite{ambrosiocolombofigalli2017, dipernalions1988b} and the references therein.

%
%
%
%
%

\subsection{Renormalized solutions}
We first recall the classical definition of renormalized solutions for transport equations.

\begin{definition} [Renormalized solution]\label{def: renormalized}
Let $\bfb :(0,T)\times \T^d \times \R^d \to \R^{2d}$ be a Borel vector field satisfying $\text{div}_{x,v}(\bfb_t) = 0$ a.e. $t\in(0,T)$ in the sense of distributions. A function $f_t$ is called a {\it renormalized solution} to the transport equation
\[
\pa_t f_t + \text{div}_{x,v}(\bfb_t f_t) = 0
\]
if, for every $\beta \in C^1(\R)\cap L^\infty(\R)$,
\[
\pa_t \beta(f_t) + \text{div}_{x,v}(\bfb_t \beta(f_t)) = 0
\]
holds in the sense of distributions. 
\end{definition}

%
%
%
%
%
\subsection{Maximal regular flows and Lagrangian solutions}
This subsection develops the concept of maximal regular flows and their connection to Lagrangian solutions.

\subsubsection{Regular flows}
By $\calL^{2d}$ we denote the Lebesgue measure in $\T^d\times\R^d$.

\begin{definition}[Regular flow]
Fix $\tau_1<\tau_2$ and $B \subseteq \T^d\times\R^d$ a Borel set. For a Borel vector field $\textbf{b}:(\tau_1,\tau_2)\times\T^d\times\R^d\to \R^{2d}$,
we say that $\bfX:[\tau_1,\tau_2]\times B\to \T^{d}\times \R^d$ is a regular flow with vector $\textbf{b}$ when
\begin{enumerate}[(i)]
\item  For a.e. $z=(x,v)\in B$, we have that $\bfX(\cdot,z)\in \AC([\tau_1,\tau_2];\T^d\times\R^d)$ and that it solves the equation $\dot{x}(t)=\textbf{b}_t(x(t))$ a.e. in $(\tau_1,\tau_2)$ with initial condition $\bfX(\tau_1,z)=z$;
\item There exists $C(\bfX) > 0$ such that $\bfX(t,\cdot)_\#(\mathcal{L}^{2d}\lefthalfcup B)\leq C\mathcal{L}^{2d}$ for all $t\in [\tau_1,\tau_2]$.
\end{enumerate}
\end{definition}

Maximal regular flows extend the idea of regular flows by accounting for particle trajectories that may escape the domain or blow up in finite time. This generalization is critical for analying systems where such phenomena cannot be ignored.

\begin{definition}[Maximal regular flow]
Let $T>0$ be some time horizon and $\bfb:(0,T)\times\T^d\times\R^d \to \R^{2d}$ a Borel vector field. For each $s\in(0,T)$, a Borel map $\bfX(\cdot,s,\cdot)$ is said to be a maximal regular flow (starting at $s$) if there exist two Borel maps $T^+_{s,\bfX}:\T^d\times\R^d\to (s,T]$, $T^-_{s,\bfX}:\T^d\times\R^d\to [0,s)$ such that $\bfX(\cdot,s,z)$ is defined in $(T^-_{s,\bfX},T^+_{s,\bfX})$ and
\begin{enumerate}[(i)]
\item For a.e. $z=(x,v)\in\T^d\times\R^d$, we have that $\bfX(\cdot,s,z)\in \AC((T^-_{s,\bfX},T^+_{s,\bfX});\T^d\times\R^d)$ and that it solves the equation $\dot{x}(t)=\textbf{b}_t(x(t))$ a.e. in $(T^-_{s,\bfX},T^+_{s,\bfX})$ with $\bfX(s,s,z)=z$.
\item There exists a constant $C = C(\bfX,s) > 0$ such that 
\[
\bfX(t,s,\cdot)_\#(\mathcal{L}^{2d}\lefthalfcup \{T^-_{s,\bfX}<t<T^+_{s,\bfX}\})\leq C\mathcal{L}^{2d}
\] 
for all $t\in [0,T]$.
\item For a.e. $z\in\T^d\times\R^d$, it either holds that 
\[
T^+_{s,\bfX}=T \mbox{ and } \bfX(\cdot,s,z)\in C([s,T];\T^d\times\R^d),\quad \mbox{or} \quad \lim_{t\uparrow T^+_{s,\bfX}}|\bfX(t,s,z)|=\infty. 
\]Analogously, either $T^-_{s,\bfX}=0$ and $\bfX(\cdot,s,z)\in C([0,s];\T^d\times\R^d)$, or $\lim_{t\downarrow T^-_{s,\bfX}}|\bfX(t,s,z)|=\infty$.
\end{enumerate}
\end{definition}

We remind the reader that existence and uniqueness of maximal regular flows were discussed in \cite{ambrosiocolombofigalli2015, ambrosiocolombofigalli2017} under varying assumptions.

The following proposition summarizes two theorems from \cite{ambrosiocolombofigalli2017}, and is also related to the ideas of \cite{BBC16}. Although originally stated for the whole space $\R^d\times\R^d$ (or an open subset thereof), the results can be naturally adapted to our setting in $\T^d\times\R^d$.

\begin{proposition}[Theorem 4.3, Theorem 4.4, \cite{ambrosiocolombofigalli2017}] \label{prop: acf}
    Let $\bfb:(0,T)\times\T^d\times\R^d \to \R^{2d}$ be a divergence-free Borel vector field satisfying the assumptions:
    \begin{enumerate}
        \item[(A1)] $\bfb \in L^1(0,T;L^1_{loc}(\T^d\times\R^d))$, 
        \item[(A2)] For any nonnegative $\bar h\in L^\infty_+(\T^d\times\R^d)$ and any closed interval $[a,b] \subset [0,T]$, both continuity equations
        \begin{align*}
            \p_t h_t \pm \textnormal{div}_{x,v} (\bfb_t h_t) = 0 \quad \text{in} \quad  (a,b) \times \T^d\times\R^d
        \end{align*}
        have at most one solution in the space of weakly$^*$ nonnegative continuous maps $[a,b]\ni t\mapsto h_t$ which verify the initial condition $h_a = \bar h$ and $\cup_{t \in [a,b]} \text{supp}\, h_t \Subset \T^d\times \R^d$. 
    \end{enumerate}
    Then there exists (up to a set of Lebesgue measure zero) a unique maximal regular flow starting at each $s\in [0,T]$ which corresponds to the vector field $\bfb$. In particular, if $\bfb_t(x,v) = (\bfb_{1t}(v),\bfb_{2t}(x))$ where
    \begin{align*}
        \bfb_1\in L^\infty(0,T;W_{loc}^{1,\infty}(\R^d;\R^d)), \qquad \bfb_{2} = \nabla G *_{x} u
    \end{align*}
    with $u\in L^\infty(0,T;\calM(\T^d))$ and $-\Delta G = \delta_0 - 1$, then $\bfb$ satisfies assumptions (A1)--(A2).
\end{proposition}
\begin{remark}
    The Green's function $G$ for the Laplacian on the torus (i.e., the function satisfying $-\Delta G=\delta_0 -1$) is of class $C^\infty(\T^d\setminus\{0\})$, and near its singularity it has the representation $G=G_1+G_2$, where $\nabla G_1= c_d x/|x|^d$ and $G_2\in C^\infty(\overline{B_{\frac{1}{4}}(0)})$ (see \cite{titchmarsh1962}). In the original setting of \cite{ambrosiocolombofigalli2017}, which considers the whole spatial space $\R^d$, the results are formulated with $\bfb_{2t} = \nabla G_1 * u_t$. However, by following the arguments of \cite[Theorem 4.4]{ambrosiocolombofigalli2017} and utilizing the Lipschitz continuity of $G_2$, it is relatively straightforward to verify that Proposition \ref{prop: acf} holds in the periodic spatial domain case.
\end{remark}

We point out that Theorem \ref{thm:pb} provides sufficient information to deduce the existence and uniqueness of a maximal regular flow corresponding to the continuity equation $\eqref{iVP}_1$, if we are provided that the bounds for the initial data in \eqref{eq:inc} propagate over time.

\begin{lemma}\label{lem: MRF unique}
Let $\textbf{b}:(0,T)\times \T^d\times\R^d\to \R^{2d}$ be given by $\textbf{b}_t(x,v)=(v,- \nabla_x \Phi[ \rho_t](x))$, where
\begin{align*}
    \rho_t = \int_{\R^d} f_t\,\dv \quad \text{and} \quad f \in L^\infty(0,T;L^1_2 \cap L^q(\T^d\times\R^d)).
\end{align*}
Then, the results of Theorem \ref{thm:pb} imply that assumption (A2) is satisfied. In particular, for any $s\in [0,T]$ there exists a unique maximal regular flow with vector $\bfb$ which starts from $s$.
\end{lemma}
\begin{proof}
Let $p := \frac{q(d-1)+2q}{q(d-1)+2}$, then we have $\rho \in L^\infty(0,T;L^p(\T^d))$. Note that Theorem \ref{thm:pb} guarantees the existence of $\Phi[\rho]\in L^\infty(0,T;W^{1,2}(\T^d))$. Due to the uniqueness part in Theorem \ref{thm:pb}, the vector field $\bfb_{2t}$ can be represented as
\[
\bfb_{2t}(x)= \nabla_x G \ast (\rho_t - e^{\Phi[\rho_t]}),  \quad\text{where}\quad  -\Delta G = \delta_0 - 1.
\]
Hence, in view of Proposition \ref{prop: acf}, it is enough to verify $\rho - e^{\Phi[\rho]}\in L^\infty(0,T;\calM(\T^d))$. Since our assumption on $\rho$ and Theorem \ref{thm:pb} provide
\[
\rho, \,\, e^{\Phi[\rho]} \in L^\infty(0,T;L^{p}(\T^d)),
\]
this is immediate.
\end{proof}

\subsubsection{Generalized flows and Lagrangian solutions}
A generalized flow is a measure concentrated on particle trajectories that evolve according to a vector field $\bfb$. It accounts for trajectories that may escape to infinity in finite time, or conversely appear from infinity and reenter the flow. Since our spatial domain is $\T^d$, note that the only way for trajectories to blow up is for their speed to become infinite. Thus, to encode the information of the trajectories that blow up, we need only introduce the one-point compactification of the velocity domain $\R^d$.

\begin{definition}[Generalized flow] \label{def: generalized flow}
    Let $\bfb:(0,T)\times\T^d\times\R^d \to \R^{2d}$ be a Borel vector field. We say $\boldsymbol{\eta}\in \calM_+(C([0,T];\T^d\times\ring{\R}^d))$ is a generalized flow of $\bfb$ if $\boldsymbol{\eta}$ is concentrated on
    \begin{align*}
        \Gamma &:= \{\eta\in C([0,T];\T^d\times\ring\R^d) : \, \eta \in \AC_{loc}(\{t:\eta(t)\notin \T^d\times\{\infty\} \} ; \T^d\times\R^d) \text{ and }\\
        & \qquad \qquad \dot\eta(t) = \bfb_t(\eta(t)) \text{ for a.e. $t$ such that $\eta(t) \notin \T^d\times\{\infty\}$}\}.
    \end{align*}
Further, let $e_t : C([0,T];\T^d\times\ring\R^d) \to \T^d\times\ring\R^d$ denote the evaluation map $e_t(\eta) = \eta(t)$. If there exists a compressibility constant $C= C(\boldsymbol{\eta})> 0$ such that
    \begin{align*}
        (e_t)\#{\boldsymbol{\eta}} \lefthalfcup (\T^d\times\R^d) \le C \calL^{2d} \quad \forall t\in [0,T],
    \end{align*}
we say $\boldsymbol{\eta}$ is regular.
\end{definition}
The one-point compactification of the velocity domain $\ring \R^d = \R^d \cup \{\infty\}$ facilitates the mathematical treatment of trajectories that escape by blowing up in velocity.

To link generalized flows with particle densities evolving under a vector field, we next define transported measures. These measures describe how densities move along the trajectories of a given maximal regular flow.

\begin{definition}[Transported measures]
Let $\boldsymbol{b}:(0,T)\times\T^d\times\R^d\to \mathbb{R}^{2d}$ be a Borel vector field having a maximal regular flow $\bfX$, and $\boldsymbol{\eta}\in\mathcal{M}_+(C[0,T];\T^d\times\ring{\R}^d)$. We say that $\boldsymbol{\eta}$ is transported by $\bfX$ if
\begin{enumerate}[(i)]
\item It holds that $(e_t)\#\boldsymbol{\eta}\ll \calL^{2d}$ for all $t\in [0,T]$.
\item For each $s\in [0,T],$ the measure $\boldsymbol{\eta}$ is concentrated on
\begin{align*}
&\{\eta\in C([0,T];\T^d\times\ring{\R}^d):\eta(s)\in \T^d\times\{\infty\} \text{ or } \eta(\cdot)=\bfX(\cdot,s,\eta(s)) \text{ in } (T^-_{s,\bfX}(\eta(s)),T^+_{s,\bfX}(\eta(s)))\}.
\end{align*}
\end{enumerate}
\end{definition}

Lagrangian solutions describe how particle densities evolve along trajectories governed by $\bfb$. This notion is foundational for linking the Eulerian and trajectory-based descriptions of particle dynamics.

\begin{definition}[Lagrangian solutions] \label{def:lag}
Assume $f \in L^\infty((0,T);L^1_{+}(\T^d\times\R^d))$ is weakly continuous on $[0,T]$ in duality with $C_c(\T^d\times\R^d)$. Let $\bfb$ be a Borel vector field admitting a maximal regular flow $\bfX$. We say that $f_t$ is a Lagrangian solution to the transport equation
\begin{align*}
    \p_t f_t + \bfb_t \cdot \nabla_{x,v} f_t = 0
\end{align*}
if there exists $\boldsymbol{\eta}\in\mathcal{M}_+(C([0,T];\T^d\times\ring{\R}^d))$ transported by $\bfX$ with $(e_t)\#\boldsymbol{\eta}=f_t\mathcal{L}^{2d}$ for every $t\in[0,T]$.
\end{definition}

\begin{remark}
    In \cite{ambrosiocolombofigalli2017}, the domain is $[0,T]\times \R^{2d}$, and thus the set above is replaced by
    \begin{align*}
        \{\eta\in C([0,T];\ring\R^{2d} : \, \eta(s) = \infty \text{ or }\eta(\cdot) = \bfX(\cdot,s,\eta(s)) \text{ in }(T_{s,\bfX}^-(\eta(s)),T_{s,\bfX}^+(\eta(s))) \}.
    \end{align*}
\end{remark}

As remarked in \cite{ambrosiocolombofigalli2017}, in case $\bfb$ is divergence-free, it is easy to observe that if $\boldsymbol{\eta}$ is transported by a maximal regular flow, then $\boldsymbol{\eta}$ is a generalized flow. On the other hand, the converse statement deserves some investigation. For the case where the spatial domain is posed as $\R^d$, it is demonstrated in \cite[Theorem 4.7]{ambrosiocolombofigalli2017}, under the assumptions that $\bfb$ is divergence-free and satisfies (A1)--(A2). Following step-by-step their procedure, it is easy to see that the results in \cite{ambrosiocolombofigalli2017} carry over precisely to our setting, and henceforth we state the proposition below without proof.

\begin{proposition}[Regular generalized flows are transported by the maximal regular flow] \label{prop: reg transp}
    Assume that $\bfb:(0,T)\times\T^d\times\R^d \to \R^{2d}$ is a divergence-free Borel vector field satisfying the assumptions (A1) and (A2) of Proposition \ref{prop: acf}. Denote by $\bfX$ its (unique) maximal regular flow. Suppose $\boldsymbol{\eta}\in \calP(C([0,T];\T^d\times\ring{\R}^d))$ is a generalized flow with respect to $\bfb$ which is regular. For $s\in [0,T]$, let $\{\boldsymbol{\eta}_z^s\} \subset \calP(C([0,T];\T^d\times\ring\R^d))$, $z\in \T^d\times\ring\R^d$, denote the disintegration of $\boldsymbol{\eta}$ with respect to the fibers of the map $e_s$ (so that $\boldsymbol{\eta}(B) = \iint_{\T^d\times\ring{\R}^d} \boldsymbol{\eta}^s_z(B) \,\textnormal{d}[(e_s)_{\#\boldsymbol{\eta}}](z)$ for all $B\subset C([0,T];\T^d\times \ring\R^d)$ Borel). Then, for $(e_s){\#\boldsymbol{\eta}}$-a.e. $z\in \T^d\times\R^d$, the measure $\boldsymbol{\eta}_z^s$ is concentrated on
    \begin{align*}
        \widehat{\Gamma}_s := \{\eta\in C([0,T];\T^d\times\ring\R^d): \eta(s) = z\text{ and } \eta(\cdot) = \bfX(\cdot,s,\eta(s)) \text{ in }(T_{s,\bfX}^-(\eta(s)), T_{s,\bfX}^+(\eta(s))) \}.
    \end{align*}
    In particular, this means $\boldsymbol{\eta}$ is transported by $\bfX$.
\end{proposition}

To represent the solution to the transport equation through generalized flow, we would like to apply the so-called \emph{superposition principle} (see \cite[Theorem 12]{ambrosiocolombofigalli2015} and \cite[Theorem 2.1]{ambrosiocrippa2008}). Note that for a direct application, global integrability conditions on the vector fields must be imposed. In the absence of such global bounds, this difficulty was overcome in \cite{ambrosiocolombofigalli2017} by introducing a ``damped" stereographic projection, as given below in Lemma \ref{lem: psi}. Since the spatial domain considered in this paper is $\T^d$, such a projection is only applicable to the velocity variable $v$. Thus, slight modification of the arguments in \cite{ambrosiocolombofigalli2017} is required to obtain an analogue, and we present the adapted results below in the fullest in order to avoid gaps. The following proposition is the analogue of the \emph{extended superposition principle} \cite[Theorem 5.1]{ambrosiocolombofigalli2017}, tailored to our spatial domain $\T^d$.

\begin{proposition}[Extended Superposition Principle on $\T^d \times \R^d$] \label{prop: re are lag}
    Let $\bfb_t(x,v) = (A(v),B_t(x))$, with $A\in L^1_{loc}(\R^d;\R^d)$ and $B\in L^1((0,T)\times\T^d;\R^d)$, be a (divergence free) Borel vector field. To the continuity equation
    \begin{align}\label{eq: transp}
        \p_t u_t + \bfb_t \cdot \nabla_{x,v} u_t = 0,
    \end{align}
    assume that $u_t \in L^\infty(0,T;$ $L^1_+(\T^d\times\R^d))$, weakly continuous on $[0,T]$ in duality with $C_c(\T^d\times\R^d)$, is either
    \begin{enumerate}[(i)]
    \item a bounded solution to \eqref{eq: transp} in the sense of distributions; or
    \item a renormalized solution to \eqref{eq: transp}.
    \end{enumerate}
    Assume also that 
    \[
    A(v)u_t\in L^1((0,T)\times\T^d\times\R^d;\R^d).
    \]
     Then there exists $\boldsymbol{\eta}\in \calM_+(C([0,T];\T^d\times\ring\R^d))$ that is concentrated on $\Gamma$ (see Definition \ref{def: generalized flow}), and which satisfies
    \begin{align*}
        &|\boldsymbol{\eta}|(C([0,T];\T^d\times\ring\R^d)) \le \sup_{t\in [0,T]} \|u_t\|_{L^1(\T^d\times\R^d)} ,\\
        &(e_t)\#\boldsymbol{\eta}\lefthalfcup(\T^d\times\R^d) = u_t \calL^{2d} \quad \forall t\in [0,T].
    \end{align*}
    Moreover, if $\bfb$ satisfies (A2) (note that (A1) is already satisfied), then $\boldsymbol{\eta}$ is transported by the maximal regular flow $\bfX$ corresponding to $\bfb$. In particular $u_t$ is a Lagrangian solution to \eqref{eq: transp}.
\end{proposition}

Before presenting the proof of this proposition, we recall the damped stereographic projection constructed in \cite{ambrosiocolombofigalli2017}. It is employed to bypass the \emph{local} integrability of $\bfb_t u_t$.

\begin{lemma}[Lemma 5.3, \cite{ambrosiocolombofigalli2017}] \label{lem: psi}
    Let $D:[0,\infty)\to (0,1]$ be a monotone nonincreasing function. Then there exist $r_0>0$ and a smooth diffeomorphism $\psi:\R^d\to \bbS^d\setminus\{N\} \subset \R^{d+1}$ which verifies
    \begin{align*}
        \begin{cases}
        \psi(v) \to N &\text{as }|v|\to\infty,\\
        |\nabla \psi(v)| \le D(0) & \forall v\in \R^d,\\
        |\nabla \psi(v)| \le D(|v|) & \forall v\in \R^d\setminus B_{r_0}.
        \end{cases}
    \end{align*}
\end{lemma}
\begin{proof}[Proof of Proposition \ref{prop: re are lag}]
For the sake of compactness of the presentation, in the sequel, we only consider the mass preserving case, i.e. $\|u_t\|_{L^1(\T^d \times \R^d)} = \|u_0\|_{L^1(\T^d \times \R^d)}$ for all $t \in [0,T]$. The general case $\|u_t\|_{L^1(\T^d \times \R^d)} \le \|u_0\|_{L^1(\T^d \times \R^d)}$ can be carried out by employing a time-dependent Dirac measure which stores the ``lost mass'', in exactly the same manner as with the proof of \cite[Theorem 5.1]{ambrosiocolombofigalli2017}.

$\bullet$ Case (i): We first consider the case where $u$ is a bounded solution to \eqref{eq: transp} in the sense of distributions. By normalizing $u_0$ we may as well assume that $\|u_t\|_{L^1(\T^d\times\R^d)} = 1$.
 
    Note that although $u_t$ is bounded, we merely have $\bfb_t u_t \in L^1_{loc}((0,T)\times \T^d\times \R^d)$, and thus no \emph{global} integrability. Toward this end, we consider $\psi$ as constructed in Lemma \ref{lem: psi}, and choose $D:[0,\infty) \to (0,1]$ such that $\psi(v)$ ensures
    \bq\label{eq:gbb}
    \int_0^T \iint_{\T^d\times\R^d} |\nabla \psi(v)| |B_t(x)| u_t(x,v)\,\dx\dv\dt < \infty.
    \eq
    To obtain this integrability, we take
    \begin{align*}
        D(r) = \begin{cases}
            1 & r\in [0,1),\\[6pt]
            \displaystyle (2^n C_n)^{-1} & r\in [2^{n-1},2^n),
        \end{cases} \quad C_n := 1+\int_0^T \int_{\T^d} \int_{B_{2^n}} |B_t(x)| u_t(x,v)\,\dv\dx\dt,
    \end{align*}
    where we note that the local integrability of $\bfb_t u_t$ provides $C_n<+\infty$. Then, by Lemma \ref{lem: psi}, we have
    \begin{align*}
        &|\nabla \psi(v)|\le 1 \quad \forall v\in \R^d,\\
        &|\nabla \psi(v)| \le (2^n C_n)^{-1} \quad \forall v\in B_{2^n}\setminus B_{2^{n-1}}, \quad \forall n \ge n_*
    \end{align*}
    for some $n_*\in \bbN$. This ensures \eqref{eq:gbb}:
        \begin{align*}
        &\int_0^T \iint_{\T^d\times\R^d} |\nabla \psi(v)| |B_t(x)| u_t(x,v)\,\dx\dv\dt \\
        &\le \int_0^T \int_{\T^d} \int_{B_{2^{n_*}}} |B_t(x)| u_t(x,v)\,\dv\dx \dt + \sum_{n=n_* + 1}^\infty 2^{-n} C_n^{-1} \int_0^T\int_{\T^d} \int_{B_{2^n\setminus B_{2^{n-1}}}} |B_t(x)| u_t(x,v)\,\dv\dx\dt \\
        &\le \int_0^T \int_{\T^d} \int_{B_{2^{n_*}}} |B_t(x)| u_t(x,v)\,\dv\dx\dt + \sum_{n=n_* + 1}^\infty 2^{-n} < +\infty.
    \end{align*}
      
    We next identify $\T^d\times \ring\R^d$ as $\T^d\times \bbS^d \subset \T^d\times \R^{d+1}$ via the diffeomorphism $\psi:\R^d\to \bbS^d\setminus \{N\}$ (namely, its extension to $\ring\R^d$), and consider the pushforward of $u_t$ to $\T^d \times \bbS^d$. We show that the corresponding pushforward measure is a solution to a continuity equation on $\T^d\times \R^{d+1}$, with corresponding vector field $\bfc_t$, which is defined below. This vector field satisfies the integrability assumptions imposed in the classical superposition principles, and therefore we apply the results of \cite{ambrosiocolombofigalli2015} to show the existence of a maximal regular flow. 
    
    Let us define $\Psi:\T^d\times\R^d \to \T^d\times (\bbS^d\setminus\{N\})$ with
\[
\Psi(x,v) = (x,\psi(v)),
\]
    and denote
    \begin{align*}
        &\mu_t := \begin{cases}
        \Psi_{\#}(u_t \calL^{2d}) & \text{on}\quad \T^d\times(\bbS^d \setminus \{N\}),\\
        0 & \text{on}\quad \T^d\times\{N\}.
        \end{cases}
    \end{align*}
    In this way, we have $\mu_t\in \calP(\T^d\times\bbS^d)$. The corresponding vector $\bfc_t$ is defined as
    \begin{align*}
        \bfc_t(y,w) := \begin{cases}
             (A(\psi^{-1}(w)), \nabla \psi(\psi^{-1}(w)) B_t(y))& \text{if }(y,w)\in \T^d\times (\bbS^d\setminus \{N\}),\\
            0 & \text{if }(y,w)\in \T^d\times \{N\}.
        \end{cases}
    \end{align*}
    We can readily confirm that $\bfc_t$ is integrable:
    \begin{align*}
    \int_0^T \iint_{\T^d\times\bbS^d} |\bfc_t| \, \textnormal{d}\mu_t \,\dt &= \int_0^T \iint_{\T^d\times \R^d} \left|(A(v), \nabla \psi(v) B_t(x) ) \right| u_t(x,v)\,\dx\dv\dt  \\
    &\le \int_0^T \iint_{\T^d\times\R^d} |A(v)| u_t(x,v)\,\dx\dv\dt  + \int_0^T \iint_{\T^d\times\R^d} |\nabla \psi(v)| |B_t(x)| u_t(x,v) \,\dx\dv\dt \\
    &< + \infty,   
    \end{align*}
    due to the assumption $Au \in L^1([0,T]\times\T^d\times\R^d)$ and \eqref{eq:gbb}.

    Then, let us show that $\mu_t$ is a distributional solution to the continuity equation corresponding to $\bfc_t$. Since $\mu_t(\T^d\times\{N\}) = 0$, for any $\varphi\in C^\infty(\T^d\times\R^{d+1})$ we have that
    \begin{align*}
        \iint_{\T^d\times\bbS^d} \varphi(y,w)\,\dd\mu_t &= \iint_{\T^d\times (\bbS^d\setminus\{N\})} \varphi(y,w) \, \dd\mu_t = \iint_{\T^d\times \R^d} \varphi(x,\psi(v)) u_t(x,v) \,\dx\dv.
     \end{align*}
     From here, we use the transport equation for $u_t$ to deduce
     \begin{align*}
         \ddt \iint_{\T^d\times \bbS^d} \varphi(y,w)\,\dd\mu_t &= \iint_{\T^d\times\R^d} \nabla_{x,v} \Big[\varphi(x,\psi(v))\Big] \cdot \bfb_t u_t \, \dx\dv \\
         &= \iint_{\T^d\times\R^d} \Big[\nabla_y \varphi(x,\psi(v)) A(v) + \nabla_w \varphi(x,\psi(v)) \nabla \psi(v) B_t(x) \Big] u_t\,\dx\dv \\
         &= \iint_{\T^d\times \bbS^d} \nabla_{y,w}\varphi(y,w) \cdot \bfc_t(y,w) \,\dd\mu_t,
     \end{align*}
     where for the final line, we use the fact that $\mu_t(\T^d\times \{N\}) = 0$. This shows that
     \begin{align*}
         \p_t \mu_t + \textnormal{div}_{y,w}(\bfc_t \mu_t) = 0.
     \end{align*}
     Consequently, applying the classical superposition principle \cite[Theorem 2.1]{ambrosiocolombofigalli2015}, we find that there exists $\boldsymbol{\sigma}\in \calP(C([0,T];\T^d\times \bbS^d))$ that is concentrated on integral curves of $\bfc$, satisfying
     \bq \label{eq: sigma}
     \mu_t = (e_t){\#\boldsymbol{\sigma}}, \quad \forall\, t\in [0,T].
     \eq

     It now only remains to pull back to $\T^d\times\ring\R^d$. Let $\tilde{\psi}:\ring\R^d \to \bbS^d$ denote the extension of $\psi$ by setting $\psi(\infty) = N$. Then, denoting $\widetilde{\Psi} = (x,\widetilde{\psi}(v))$, we define the map $\Xi:C([0,T];\T^d\times\bbS^d) \to C([0,T];\T^d\times \ring\R^d)$ by
     \begin{align*}
         \Xi(\nu)(t) = \Big(\widetilde{\Psi}^{-1} \circ \nu\Big)(t).
     \end{align*}
     The measure
     \begin{align*}
         \boldsymbol{\eta} := \Xi{\#\boldsymbol{\sigma}}
     \end{align*}
    is then concentrated on $\Gamma$ of Definition \ref{def: generalized flow}, and is therefore a generalized flow of $\bfb$, which satisfies
    \begin{equation*}
        |\boldsymbol{\eta}|(C([0,T];\T^d\times\ring\R^d)) \le 1 = \|u_t\|_{L^1(\T^d\times\R^d)}.
    \end{equation*}
    Finally, we observe from \eqref{eq: sigma} that
    \begin{equation} \label{eq: ut}
    \begin{split}
        (e_t){\#\boldsymbol{\eta}}\lefthalfcup (\T^d\times\R^d) &= (e_t){\# \Xi \# \boldsymbol{\sigma}} \lefthalfcup(\T^d\times\R^d) \\
        &= \widetilde{\Psi}^{-1} \# e_t \# \boldsymbol{\sigma} \lefthalfcup (\T^d\times\R^d) \\
        &= \widetilde{\Psi}^{-1}\# \mu_t \lefthalfcup(\T^d\times\R^d)\\
        &= u_t \calL^{2d}.
    \end{split}
    \end{equation}
    In particular, we have
    \begin{align*}
        (e_t)\#\boldsymbol{\eta} \lefthalfcup(\T^d\times\R^d) \le \|u\|_{L^\infty} \calL^{2d},
    \end{align*}
    and therefore $\boldsymbol{\eta}$ is a regular generalized flow. If $\bfb$ satisfies (A2), then its maximal regular flow $\bfX$ is unique, and by Proposition \ref{prop: reg transp} $\boldsymbol{\eta}$ is transported by $\bfX$. By definition, we note that \eqref{eq: ut} implies $u_t$ is a Lagrangian solution to the transport equation, as desired.

$\bullet$ Case (ii): In the case of renormalized solutions, we note that for any $\beta\in C^1\cap L^\infty(\R)$, $\beta(u_t)$ is a bounded solution to the transport equation \eqref{eq: transp} in the sense of distributions. Hence, the above results apply to each $\beta(u_t)$. By approximation arguments, one readily checks that the renormalization property also holds true for bounded and Lipschitz $\beta$. Thus, by setting
\begin{align*}
    \beta_k(u) = \begin{cases}
        0 & u\le k,\\
        u-k & k \le u \le k+1,\\
        1 & u\ge k+1.
    \end{cases}
\end{align*}
By case (i), for each $k$, there is a measure $\boldsymbol{\eta}_k\in \calM_+(C([0,T];\T^d\times\ring\R^d))$ which is concentrated on the set defined in Definition \ref{def: generalized flow}, and which satisfies
\begin{align*}
    &|\boldsymbol{\eta}_k|(C([0,T];\T^d\times\ring\R^d)) \le \sup_{t\in [0,T]} \|\beta_k(u_t)\|_{L^1(\T^d\times\R^d)},\\
    &(e_t)\#\boldsymbol{\eta}_k \lefthalfcup (\T^d\times\R^d) = \beta_k(u_t) \calL^{2d} \quad \forall t\in [0,T].
\end{align*}
By observing that
\begin{align*}
    \sum_{k=0}^\infty \beta_k(u_t) = u_t,
\end{align*}
it is immediate that the measure $\boldsymbol{\eta} := \sum_{k=0}^\infty \boldsymbol{\eta}_k$ satisfies all the assertions of the proposition. This concludes the proof.
\end{proof}


%
\bibliographystyle{abbrv}

\end{document}